\documentclass[10pt]{article}
\usepackage{mathrsfs}
\usepackage{amsthm}
\usepackage{amssymb}
\usepackage{amsmath}
\usepackage{graphicx}
\usepackage{color}
\usepackage{amsfonts}
\usepackage{float}
\usepackage{cite}
\usepackage[text={140mm,210mm},left=45mm,vmarginratio=1:1]{geometry}
\newtheorem{theorem}{Theorem}[section]

\newtheorem{lemma}[theorem]{Lemma}
\newtheorem{proposition}[theorem]{Proposition}

\numberwithin{equation}{section}
\normalsize

\begin{document}
\title{\textbf{Moderate deviations of density-dependent Markov Chains}}

\author{Xiaofeng Xue \thanks{\textbf{E-mail}: xfxue@bjtu.edu.cn \textbf{Address}: School of Science, Beijing Jiaotong University, Beijing 100044, China.}\\ Beijing Jiaotong University}

\date{}
\maketitle

\noindent {\bf Abstract:} The density-dependent Markov chain (DDMC) introduced in \cite{Kurtz1978} is a continuous time Markov process applied in fields such as epidemics, chemical reactions and so on. In this paper, we give moderate deviation principles of paths of DDMC under some generally satisfied assumptions. The proofs for the lower and upper bounds of our main result utilize an exponential martingale and a generalized version of Girsanov's theorem. The exponential martingale is defined according to the generator of DDMC.

\quad

\section{Introduction}\label{section one}
In this paper we are concerned with the density-dependent Markov process (DDMC) introduced in \cite{Kurtz1978}. For each integer $n\geq 1$, the density-dependent Markov chain $\left\{X_t^n\right\}_{t\geq 0}$ is a Markov process with state space $nG$, where $G\subseteq \mathbb{R}^d$ is a given closed and convex set and
\[
nG=\left\{y\in \mathbb{R}^d:~\frac{y}{n}\in G\right\}.
\]
The transition rates function of $\{X_t^n\}_{t\geq 0}$ is given by
\[
X_t^n\rightarrow X_t^n+l \text{~at rate~} nF_l\left(\frac{X_t^n}{n}\right)
\]
for each $l\in \mathcal{A}$, where $\mathcal{A}$ is a subset of $\mathbb{R}^d$ while $F_l\in C^1(\mathbb{R}^d)$ for each $l\in \mathcal{A}$. To ensure
$P(X_t^n\in nG\text{~for all~}t)=1$, $\{F_l\}_{l\in \mathcal{A}}$ further satisfy
\[
F_l(x)=0
\]
if $nx\in nG$ but $nx+l\not\in nG$ for some $n\geq 1$.

Important examples of DDMC are given in former references, such as \cite{Chan1998, Ethier1986, Ge2017, Kurtz1978, Pardoux2017} and so on. Here we recall some of these examples. Note that we consider elements of $\mathbb{R}^d$ as column vectors for later use while we use $\mathsf{T}$ to denote the transposition operator.

\textbf{Example 1} \emph{The contact process on the complete graph}. Let $d=1$, $\lambda>0$, $G=[0,1]$, $\mathcal{A}=\{1,-1\}$  and
\[
F_1(x)=\lambda x(1-x), \text{~} F_{-1}(x)=x
\]
for $x\in G$ (we do not care $F_{\pm1}(x)$ for $x\not\in G$), then $X_t^n$ is the number of infected vertices at moment $t$ for the contact process with infection rate $\lambda/n$ on the complete graph with $n$ vertices.

For the contact process on the complete graph, each vertex is healthy or infected. An infected vertex recovers at rate one while a healthy vertex is infected at rate proportional to the number of infected vertices. For a detailed survey of the study of the contact process, see Chapter 6 of \cite{Lig1985}.

\qed

\quad

\textbf{Example 2} \emph{The SIR model on the complete graph}. Let $d=2$, $\lambda>0$,
\[
G=\{(x,y)^{\mathsf{T}}:x,y\geq 0, x+y\leq 1\}, \text{~}\mathcal{A}=\{(0,-1)^{\mathsf{T}},(-1,1)^{\mathsf{T}}\}
\]
and
\[
F_{_{(0,-1)^{\mathsf{T}}}}(x,y)=y, \text{~} F_{_{(-1,1)^{\mathsf{T}}}}(x,y)=\lambda xy
\]
for $(x,y)^{\mathsf{T}}\in G$, then $X_t^n=(S_t^n,I_t^n)^{\mathsf{T}}$ is the state at moment $t$ of the SIR model with infection rate $\lambda/n$ on the complete graph with $n$ vertices, where $S_t^n$ is the number of susceptible vertices while $I_t^n$ is the number of infected vertices.

For the SIR model, which is also called as the epidemic model, a vertex is in one of the three states `susceptible', `infected' and `removed'. An infected vertex is removed at rate one while a susceptible vertex is infected at rate proportional to the number of infected vertices.

\qed

\quad

\textbf{Example 3} \emph{Chemical reactions}. Here we only discuss a special simple case. For general cases, see Section 11.1 of \cite{Ethier1986} or \cite{Ge2017}. Assuming $3$ chemical reactants $R_1$, $R_2$, $R_3$ are undergoing the chemical reaction
\[
R_1+R_2\rightleftharpoons R_3
\]
in a system with at most $n$ molecules. If the forward reaction occurs at rate $\lambda/n$ for a given pair of a $R_1$ molecule and a $R_2$ molecule while the reverse reaction occurs at rate $\mu$ for
a given $R_3$ molecule, then this chemical reaction can be described by DDMC $\left\{X_t^n=\left(X_t^{n,1},X_t^{n,2},X_t^{n,3}\right)^{\mathsf{T}}:~t\geq 0\right\}$ with $d=3$,
\[
G=\{(x,y,z)^{\mathsf{T}}:x,y,z\geq 0, x+y+2z\leq 1\}, \text{~}\mathcal{A}=\big\{(-1,-1,1)^{\mathsf{T}}, (1,1,-1)^{\mathsf{T}}\big\}
\]
and
\[
F_{_{(-1,-1,1)^{\mathsf{T}}}}(x,y,z)=\lambda xy, \text{~} F_{_{(1,1,-1)^{\mathsf{T}}}}(x,y,z)=\mu z,
\]
where $X_t^{n,i}$ is the number of $R_i$ molecules at moment $t$.

\qed

\textbf{Example 4} \emph{Yule process}. The Yule process with rate $\lambda$ is also a DDMC with $d=1$,
\[
G=[0,+\infty), \text{~}\mathcal{A}=\big\{1\big\}
\]
and $F_1(x)=\lambda x$.

\qed

\quad

Law of large numbers (LLN) and central limit theorem (CLT) of DDMC are given in \cite{Kurtz1978}.
\begin{proposition}\label{Proposition 1.1 LLN of DDMC} (Kurtz, 1978)
If $\sum_{l\in \mathcal{A}}lF_l(x)$ satisfies Lipschitz condition on $G$ while
\[
\lim_{n\rightarrow+\infty}\frac{X_0^n}{n}=x_0 \text{~in probability},
\]
then
\[
\lim_{n\rightarrow +\infty}\left(\sup_{0\leq t\leq T_0}\left\|\frac{X_t^n}{n}-X_t\right\|\right)=0
\]
in probability for any $T_0>0$, where $\{X_t\}_{t\geq 0}$ is the unique solution to the ODE
\begin{equation}\label{equ 1.1 ODE}
\begin{cases}
&\frac{d}{dt}X_t=\sum_{l\in \mathcal{A}}lF_l(X_t),\\
&X_0=x_0
\end{cases}
\end{equation}
and $\left\|x\right\|$ is the $L^1$ norm of $x$.
\end{proposition}

To recall CLT theorem, let $Y_t^n=\frac{X_t^n-nX_t}{\sqrt{n}}$ for any $t\geq 0$.
\begin{proposition}\label{Proposition 1.2 CLT of DDMC}(Kurtz, 1978)
Under the assumption of Proposition \ref{Proposition 1.1 LLN of DDMC},
if $Y_0^n$ converges weakly to $V_0$ as $n\rightarrow+\infty$, then
$\{Y_t^n:~0\leq t\leq T_0\}$ converges weakly to $\{V_t:~0\leq t \leq T_0\}$ as $n\rightarrow+\infty$, where $\{V_t\}_{t\geq 0}$ is a time-inhomogeneous O-U process:
\[
dV_t=d\alpha_t+\sum_{l\in\mathcal{A}}l(\nabla F_l(X_t)\cdot V_t)dt
\]
such that
\[
\alpha_t=\sum_{l\in \mathcal{A}}lW_l\left(\int_0^tF_l(X_s)ds\right),
\]
where $\{W_l\}_{l\in \mathcal{A}}$ are independent standard Brownian motions and $\nabla=\left(\frac{\partial}{\partial x_1},\ldots,\frac{\partial}{\partial x_d}\right)^{\mathsf{T}}$.

\end{proposition}
Large deviations are also discussed for DDMC. Under different assumptions of $\mathcal{A}$ and $\{F_l\}_{\mathcal{A}}$, large deviations of paths of DDMC are established in Chapter 5 of \cite{Schwartz1995} and \cite{Agazzi2018, Chan1998, Pardoux2017} respectively.

In this paper we are concerned with the moderate deviation of DDMC, i.e., the goal of this paper is to show that
\[
P\left(\left\{\frac{X_t^n-nX_t}{a_n}:0\leq t\leq T_0\right\}=dx\right)\approx \exp\left(-\frac{a_n^2}{n}I(x)\right)dx
\]
under some generally satisfied assumptions for any $x:[0,T_0]\rightarrow \mathbb{R}^d$ in the Skorokhod space $\mathcal{D}\left([0,T_0], \mathbb{R}^d\right)$ and sequence $\{a_n\}_{n\geq 1}$ satisfying
\[
\frac{a_n}{n}\rightarrow 0, \text{~}\frac{a_n}{\sqrt{n}}\rightarrow+\infty
\]
with a rate function $I:\mathcal{D}\left([0,T_0], \mathbb{R}^d\right)\rightarrow [0,+\infty]$. References (see page 285 of \cite{Deuschel1989} or page 577 of \cite{Gao2003}) show that the study of the moderate deviation dates back to 1928, when Khinchin gives the moderate deviation for independent Bernoulli stochastic variables. Over several past decades, moderate deviations are obtained for many different types of stochastic processes. References \cite{Borovkov1978, Borovkov1980, deAcosta1998, Gao1996, Gao2003, Gao2010, Gao2017, Wang2016, WangR2015, Wu1995} and so on can be consulted for an outline of this development.

This paper is inspired a lot by \cite{Gao2003}, where the moderate deviation of the hydrodynamic limit of the symmetric exclusion process (SEP) is discussed. Evidences show that DDMC has limit behavior similar with the hydrodynamic of the SEP. As we have recalled, LLN of DDMC is driven by an ODE on $\mathbb{R}^d$ while the hydrodynamic of the SEP is driven by a heat equation, which can be considered as an ODE on the space of measures (see \cite{Kipnis1989}). CLT of DDMC is driven by a time-inhomogeneous O-U process on $\mathbb{R}^d$ while CLT of the hydrodynamic of the SEP is driven by a time-inhomogeneous O-U process on the space of measures (see Chapter 11 of \cite{Kipnis1999}). As a result, we are motivated to study moderate deviation of DDMC, which is expected to be an analogue of the main result given in \cite{Gao2003}.

\section{Main results}\label{section two main result}
In this section we give our main results. First we introduce some notations and basic assumptions for later use. For any $x=(x_1,\ldots,x_d)^{\mathsf{T}}\in \mathbb{R}^d$, we use $\left\|x\right\|$ to denote the $L^1$ norm of $x$, i.e., $\left\|x\right\|=\sum_{i=1}^d|x_i|$.

Throughout this paper, we adopt the following basic assumptions.

Assumption (1): $x_0\neq 0$ is a given point in $G$.

Assumption (2): for each $n\geq 1$, $X_0^n=nx_0$.

Assumption (3): $\{a_n\}_{n\geq 1}$ is a positive sequence such that $\frac{a_n}{n}\rightarrow 0$ and $\frac{a_n}{\sqrt{n}}\rightarrow +\infty$ as $n\rightarrow +\infty$.

Assumption (4): $\mathcal{A}$ is finite.

Assumption (5): $F_l(0)=0$ for each $l\in \mathcal{A}$ and there exists $K_1<+\infty$ such that $\left\|\nabla F_l(x)\right\|\leq K_1$ for any $l\in \mathcal{A}$ and $x\in G$.

\quad

It is easy to check that all the four examples in Section \ref{section one} satisfies Assumptions (1)-(5). Note that we do not assume that $F_l$ is bounded on $G$ (which Examples 1-3 satisfies) to make our results can be applied in examples such as Yule processes, where $F_l(x)$ is dominated from above by a linear function of $\left\|x\right\|$ but unbounded on $G$.

For given $T_0>0$, we use $\mathcal{D}\left([0,T_0], \mathbb{R}^d\right)$ to denote the set of c\`{a}dl\`{a}g functions $f: [0,T_0]\rightarrow \mathbb{R}^d$ with $f_0=0$. Under the metric introduced in \cite{Skorokhod1956}, $\mathcal{D}\left([0,T_0], \mathbb{R}^d\right)$ is a complete separable metric space, i.e, Skorokhod space.

Now we give the rate function. For any $t\geq 0$, we define
\[
b_t=\sum_{l\in \mathcal{A}}l(\nabla^{\mathsf{T}}F_l)(X_t) \text{\quad and \quad} \sigma_t=\sum_{l\in \mathcal{A}}lF_l(X_t)l^{\mathsf{T}},
\]
where $\{X_t\}_{t\geq 0}$ is defined as in Equation \eqref{equ 1.1 ODE} and $\nabla^{\mathsf{T}}=\left(\frac{\partial}{\partial x_1},\ldots,\frac{\partial}{\partial x_d}\right)$. Note that $b_t, \sigma_t$ are both $d\times d$ matrices.
Then, for any $f\in \mathcal{D}\left([0,T_0],\mathbb{R}^d\right)$, we define
\begin{align} \label{equ 2.1 definition of rate function}
I(f)=\sup\big\{&f(T_0)\cdot g(T_0)-\int_0^{T_0}f_s\cdot g_s^{\prime}ds-\int_0^{T_0}\left(b_sf_s\right)\cdot g_sds \notag\\
&-\frac{1}{2}\int_0^{T_0} g_s^{\mathsf{T}}\sigma_sg_sds: ~g\in C^2\left([0,T_0], \mathbb{R}^d\right)\big\},
\end{align}
where
\[
g^{\prime}(t)=\left(g_1^{\prime}(t),g_2^{\prime}(t),\ldots,g_d^{\prime}(t)\right)^{\mathsf{T}}
\]
for any $g(t)=\left(g_1(t),g_2(t),\ldots,g_d(t)\right)^{\mathsf{T}} \in C^2\left([0,T_0], \mathbb{R}^d\right)$ and $x\cdot y=\sum_{i=1}^dx_iy_i$ for $x=(x_1,\ldots, x_d)^{\mathsf{T}}$ and $y=(y_1,\ldots, y_d)^{\mathsf{T}}$.

Now we give our main result. For simplicity, we use $\vartheta^n$ to denote the path of $\big\{\frac{X_t^n-nX_t}{a_n}:0\leq t\leq T_0\big\}$.
\begin{theorem}\label{theorem main 2.1 MDP} Under Assumptions (1)-(5), for any open set $O \subseteq \mathcal{D}\left([0,T_0],\mathbb{R}^d\right)$,
\[
\liminf_{n\rightarrow+\infty}\frac{n}{a_n^2}\log P\left(\vartheta^n\in O\right)\geq -\inf_{f\in O}I(f),
\]
while for any closed set $C \subseteq \mathcal{D}\left([0,T_0], \mathbb{R}^d\right)$,
\[
\limsup_{n\rightarrow+\infty}\frac{n}{a_n^2}\log P\left(\vartheta^n\in C\right)\leq -\inf_{f\in C}I(f),
\]
where $I$ is defined as in Equation \eqref{equ 2.1 definition of rate function}. Furthermore, if $\sigma(t)$ is invertible for $0\leq t\leq T_0$, then
\begin{equation}\label{equ 2.2}
I(f)=
\begin{cases}
\frac{1}{2}\int_0^{T_0} (f_s^{\prime}-b_sf_s)^{\mathsf{T}}\sigma^{-1}_s(f_s^{\prime}-b_sf_s)ds & \text{~if~}f \text{~is absolutely continuous},\\
+\infty & \text{~otherwise}.
\end{cases}
\end{equation}

\end{theorem}

Here we give an intuitive explanation of Theorem \ref{theorem main 2.1 MDP} in the case where $d=1$ and $\sigma_t\neq 0$. By Proposition \ref{Proposition 1.2 CLT of DDMC}, $\frac{X_t^n-nX_t}{a_n}\approx \frac{\sqrt{n}}{a_n}V_t$, where $\{V_t\}_{t\geq 0}$ is the solution of
\[
\begin{cases}
& dV_t=\sqrt{\sigma_t}dW_t+b_tV_tdt,\\
& V_0=0.
\end{cases}
\]
Then, it is natural to non-rigorously think
\[
P\left(\vartheta^n=df\right)\approx P\left(V_\cdot=d\frac{a_nf}{\sqrt{n}}\right).
\]
Let $0=t_0<t_1<\ldots<t_M=T_0$ be a partition of $[0,T_0]$ with $\sup_{i}(t_{i+1}-t_i)$ very small, then $V_\cdot=d\frac{a_nf}{\sqrt{n}}$ can be non-rigorously interpreted as
\[
\frac{a_n}{\sqrt{n}}f_{t_{i+1}}-\frac{a_n}{\sqrt{n}}f_{t_{i}}\approx b_{t_i}\frac{a_n}{\sqrt{n}}f_{t_{i}}(t_{i+1}-t_i)+\sqrt{\sigma_{t_i}}(W_{t_{i+1}}-W_{t_i}),
\]
i.e.,
\[
\sqrt{\sigma_{t_i}}(W_{t_{i+1}}-W_{t_i})\approx \frac{a_n}{\sqrt{n}}(t_{i+1}-t_i)(f^\prime_{t_i}-b_{t_i}f_{t_i})
\]
for each $i$. Since $\sqrt{\sigma_{t_i}}(W_{t_{i+1}}-W_{t_i})$ follows from $N\left(0,\sigma_{t_i}(t_{i+1}-t_i)\right)$, the above event occurs with probability about
\[
\exp\big\{-\frac{a_n^2}{n}\frac{\left(f^\prime_{t_i}-b_{t_i}f_{t_i}\right)^2(t_{i+1}-t_i)^2}{2\sigma_{t_i}(t_{i+1}-t_i)}\big\}dx
=\exp\big\{-\frac{a_n^2}{n}\frac{\left(f^\prime_{t_i}-b_{t_i}f_{t_i}\right)^2(t_{i+1}-t_i)}{2\sigma_{t_i}}\big\}dx.
\]
Since $\{W_{t_{i+1}}-W_{t_i}\}_i$ are independent, $V_\cdot=d\frac{a_n}{\sqrt{n}}f$ occurs with probability about
\[
\exp\big\{-\frac{a_n^2}{n}\sum_{i}\frac{\left(f^\prime_{t_i}-b_{t_i}f_{t_i}\right)^2(t_{i+1}-t_i)}{2\sigma_{t_i}}\big\}dx
\approx \exp\big\{-\frac{a_n^2}{n}\int_0^{T_0}\frac{(f^\prime_{s}-b_{s}f_{s})^2}{2\sigma_s}ds\big\}dx,
\]
which non-rigorously shows that the rate function
\[
I(f)=\int_0^{T_0}\frac{(f^\prime_{s}-b_{s}f_{s})^2}{2\sigma_s}ds.
\]

\quad

The rigorous proofs of the lower and upper bounds in Theorem \ref{theorem main 2.1 MDP} are given in Sections \ref{section four} and \ref{section five}
respectively. The strategy of our proofs is similar with that utilized in \cite{Gao2003}, where an exponential martingale will be introduced. To define this martingale rigorously, some basic properties of $\{X_t^n\}_{0\leq t\leq T_0}$ are given in Section \ref{section three}.

At the end
of this section, we show that the two definitions of the rate function $I$ given in Equations \eqref{equ 2.1 definition of rate function} and \eqref{equ 2.2} are
equivalent under the assumption that $\{\sigma_t\}_{0\leq t\leq T_0}$ is invertible.

\proof[Proof of Equation \eqref{equ 2.2}]

For $f\in \mathcal{D}\left([0,T_0],\mathbb{R}^d\right)$, we only need to show that $I(f)<+\infty$ implies that $f$ is absolutely continuous and
\[
I(f)=\frac{1}{2}\int_0^{T_0} (f_s^{\prime}-b_sf_s)^{\mathsf{T}}\sigma^{-1}_s(f_s^{\prime}-b_sf_s)ds.
\]
For $f$ makes $I(f)<+\infty$, we define
\[
\mathcal{L}_{1,f}(g)=f(T_0)\cdot g(T_0)-\int_0^{T_0}f_s\cdot g_s^{\prime}ds-\int_0^{T_0}\left(b_sf_s\right)\cdot g_sds
\]
and
\[
\mathcal{L}_2(g)=\int_0^{T_0} g_s^{\mathsf{T}}\sigma_sg_sds
\]
for each $g\in C^2\left([0,T_0], \mathbb{R}^d\right)$. Then
\[
I(f)=\sup\big\{\mathcal{L}_{1,f}(g)-\frac{1}{2}\mathcal{L}_2(g):g\in C^2\left([0,T_0], \mathbb{R}^d\right)\big\}.
\]
For each $c\in \mathbb{R}$ and $g\neq0$,
\[
\mathcal{L}_{1,f}(cg)-\frac{1}{2}\mathcal{L}_2(cg)=c\mathcal{L}_{1,f}(g)-c^2\frac{1}{2}\mathcal{L}_2(g).
\]
Hence, $\mathcal{L}_{1,f}(cg)-\frac{1}{2}\mathcal{L}_2(cg)$ get the maximum $\frac{\left(\mathcal{L}_{1,f}(g)\right)^2}{2\mathcal{L}_2(g)}$ when
$c=\frac{\mathcal{L}_{1,f}(g)}{\mathcal{L}_2(g)}$. As a result,
\begin{equation}\label{equ 2.3}
I(f)=\sup\big\{\frac{\left(\mathcal{L}_{1,f}(g)\right)^2}{2\mathcal{L}_2(g)}:g\in C^2\left([0,T_0], \mathbb{R}^d\right) \text{~and~} g\neq 0\big\}.
\end{equation}
Let $L^2_\sigma\left([0,T],\mathbb{R}^d\right)$ be the set of measurable $g:[0,T_0]\rightarrow \mathbb{R}^d$ such that
\[
\int_0^{T_0}g_s^{\mathsf{T}}\sigma_s g_sds<+\infty.
\]
Under the assumption that $\sigma_t$ is invertible for $0\leq t\leq T_0$, $\sigma_t$ is positive-definite for $0\leq t\leq T_0$.
Therefore, $L^2_\sigma\left([0,T], \mathbb{R}^d\right)$ is a Hilbert space under the inner product
\[
\langle f,g\rangle_{L^2_\sigma}=\int_0^{T_0}f_s^{\mathsf{T}}\sigma_s g_sds
\]
for $f,g\in L^2_\sigma\left([0,T], \mathbb{R}^d\right)$. Note that
\[
\sqrt{\mathcal{L}_2(g)}=\sqrt{\langle g,g\rangle_{L^2_\sigma}},
\]
which is the norm of $g$ generated by $\langle \cdot,\cdot\rangle_{L^2_\sigma}$. For any $g\in C^2\left([0,T_0], \mathbb{R}^d\right)$, by Equation \eqref{equ 2.3},
\[
|\mathcal{L}_{1,f}(g)|\leq \sqrt{2I(f)}\sqrt{\langle g,g\rangle_{L^2_\sigma}}.
\]
Since $C^2\left([0,T_0], \mathbb{R}^d\right)$ is dense in $L^2_\sigma\left([0,T], \mathbb{R}^d\right)$, $\mathcal{L}_{1,f}$ can be extended to
\[
\widetilde{\mathcal{L}}_{1,f}: L^2_\sigma\left([0,T], \mathbb{R}^d\right)\rightarrow \mathbb{R}
\]
such that $\widetilde{\mathcal{L}}_{1,f}\Big|_{C^2\left([0,T_0], \mathbb{R}^d\right)}=\mathcal{L}_{1,f}$ and
\[
|\widetilde{\mathcal{L}}_{1,f}(g)|\leq \sqrt{2I(f)}\sqrt{\langle g,g \rangle_{L^2_\sigma}}
\]
for any $g\in L^2_\sigma\left([0,T], \mathbb{R}^d\right)$. That is to say, $\widetilde{\mathcal{L}}_{1,f}$ is a bounded linear operator on $L^2_\sigma([0,T], \mathbb{R}^d)$.
Therefore, according to Riesz representation theorem, there exists $\psi\in L^2_\sigma\left([0,T], \mathbb{R}^d\right)$ such that
\begin{equation}\label{equ 2.4}
\widetilde{\mathcal{L}}_{1,f}(g)=\langle g,\psi\rangle_{L^2_\sigma}
\end{equation}
for any $g\in L^2_\sigma\left([0,T], \mathbb{R}^d\right)$. As a result, by the definition of $\mathcal{L}_{1,f}$,
\[
f(T_0)\cdot g(T_0)-f(0)\cdot g(0)-\int_0^{T_0}f_s\cdot g_s^{\prime}ds=\int_0^{T_0}(b_sf_s+\sigma_s\psi_s)\cdot g_s ds
\]
for each $g\in C^2\left([0,T], \mathbb{R}^d\right)$. Then, according to the formula of integration by parts, $f$ is absolutely continuous and
\[
f_t^\prime=b_tf_t+\sigma_t\psi_t,
\]
i.e.,
\begin{equation}\label{equ 2.5}
\psi_t=\sigma_t^{-1}(f_t^\prime-b_tf_t).
\end{equation}
By Equation \eqref{equ 2.4} and Cauchy-Schwartz inequality, for any $g\in C^2\left([0,T], \mathbb{R}^d\right)$,
\begin{equation*}
\left(\mathcal{L}_{1,f}(g)\right)^2\leq \langle g,g \rangle_{L^2_\sigma}\langle\psi,\psi\rangle_{L^2}=\mathcal{L}_2(g)\langle\psi,\psi\rangle_{L^2_\sigma}.
\end{equation*}
Therefore, by Equation \eqref{equ 2.3},
\begin{equation}\label{equ 2.6}
I(f)\leq \frac{1}{2}\langle\psi,\psi\rangle_{L^2_\sigma}.
\end{equation}
On the other hand, let $g_n\in C^2\left([0,T], \mathbb{R}^d\right)$ such that $\lim_{n\rightarrow+\infty}g_n=\psi$ under the distance generated by $\langle \cdot,\cdot
\rangle_{L^2_\sigma}$, then by Equation \eqref{equ 2.4},
\[
\mathcal{L}_{1,f}(g_n)=\langle g_n,\psi\rangle_{L^2_\sigma}\rightarrow \langle \psi,\psi\rangle_{L^2_\sigma} \text{~and~}\langle g_n,g_n\rangle_{L^2_\sigma}\rightarrow \langle \psi,
\psi\rangle_{L^2_\sigma}
\]
while
\[
\frac{\left(\mathcal{L}_{1,f}(g_n)\right)^2}{2\mathcal{L}_2(g_n)}=\frac{\left(\mathcal{L}_{1,f}(g_n)\right)^2}{2\langle g_n,g_n\rangle_{L^2_\sigma}}\rightarrow
\frac{1}{2}\langle \psi,\psi\rangle_{L^2_\sigma}.
\]
Hence, by Equation \eqref{equ 2.3},
\begin{equation}\label{equ 2.7}
I(f)\geq \frac{1}{2}\langle \psi,\psi\rangle_{L^2_\sigma}.
\end{equation}
By Equations \eqref{equ 2.5}, \eqref{equ 2.6} and \eqref{equ 2.7},
\[
I(f)=\frac{1}{2}\langle \psi,\psi\rangle_{L^2_\sigma}=\frac{1}{2}\int_0^{T_0}(f_s^\prime-b_sf_s)^{\mathsf{T}}\sigma_s^{-1}(f_s^\prime-b_sf_s)ds.
\]

\qed

\section{Preliminary results}\label{section three}
In this section we give some preliminary results of $\{X_t^n\}_{t\geq 1}$ for later applications in the proof of Theorem \ref{theorem main 2.1 MDP}, i.e., the goal of this section is to prove the following three lemmas.

\begin{lemma}\label{lemma 3.1} Under Assumptions (1)-(5), there exists $\theta>0$ such that
\[
E\left(\exp\big\{\Lambda \sup\limits_{0\leq t\leq T_0}\left\|X_t^n\right\|\big\}\right)<+\infty
\]
for any $\Lambda\in [0,\theta)$ and all $n\geq 1$.
\end{lemma}

\begin{lemma}\label{lemma 3.2} Under Assumptions (1)-(5), there exist $K_2, K_3\in (0, +\infty)$ such that
\[
P\left(\sup_{0\leq t\leq T_0}\left\|X_t^n\right\|>nK_2\right)\leq \exp\{-K_3n\}
\]
for all $n\geq 1$.
\end{lemma}

\begin{lemma}\label{lemma 3.3}Under Assumptions (1)-(5), for any $\epsilon>0$, there exists $K_4(\epsilon)\in (0,+\infty)$ such that
\[
P\left(\sup_{0\leq t\leq T_0}\left\|\frac{X_t^n}{n}-X_t\right\|>\epsilon\right)\leq \exp\big\{-K_4(\epsilon)n\big\}
\]
for sufficiently large $n$, where $\{X_t\}_{t\geq 0}$ is defined as in Equation \eqref{equ 1.1 ODE}.
\end{lemma}

Readers may think that Lemmas \ref{lemma 3.2} and \ref{lemma 3.3} are corollaries of the large deviation principle of DDMC given in \cite{Agazzi2018, Chan1998, Pardoux2017} or \cite{Schwartz1995}. However, the main theories in \cite{Agazzi2018, Chan1998, Pardoux2017} or \cite{Schwartz1995} requires the assumption that $\{F_l(x)\}_{l\in \mathcal{A}}$ are bounded on $\mathbb{R}^d$. Hence the proofs of Lemmas \ref{lemma 3.2} and \ref{lemma 3.3} are still needed under our assumptions (1)-(5).

Note that Lemmas \ref{lemma 3.1}-\ref{lemma 3.3} relies heavily on the assumption that $\mathcal{A}$ is finite. Estimations of moments of $X_t^n$ under a general assumption where $\mathcal{A}$ is infinite can be found in Theorems 2.1 and 2.2 of \cite{Kurtz1978}.

\proof[Proof of Lemma \ref{lemma 3.1}] Since $\mathcal{A}$ is finite, by Assumption (5), there exist $K_6, K_7\in (0,+\infty)$ such that $\{\left\|X_t^n\right\|\}_{t\geq 0}$ is stochastically dominated from above by the Markov process $\{\eta_t^n\}_{t\geq 0}$ with $\eta_0^n=\left\|X_0^n\right\|$ and transition rates function given by
\[
\eta_t^n\rightarrow \eta_t^n+K_7 \text{~at rate~}K_6\eta_t^n.
\]
Without loss of generality, we assume that $x_0(i)/K_7$ is an integer for $1\leq i\leq d$.
For each $n\geq 1$, we use $\{\widetilde{\eta}_t^n\}_{t\geq 0}$ to denote the Yule process with rate $1$ and initial state $\widetilde{\eta}_0^n=n$, i.e.,
\[
\widetilde{\eta}_t^n\rightarrow \widetilde{\eta}_t^n+1 \text{~at rate~} \widetilde{\eta}_t^n.
\]
Then,
$
\big\{\frac{1}{K_7}\eta^n_{_{t/(K_6K_7)}}:t\geq 0\big\}
$
is a copy of $\big\{\widetilde{\eta}_t^{n\left\|x_0\right\|/K_7}\big\}_{t\geq 0}$. By classic theory of Yule process, $\widetilde{\eta}^1_t$ follows geometric distribution with parameter $e^{-t}$ while $\widetilde{\eta}^n_{t}$ can be written as
\begin{equation}\label{equ 3.1}
\widetilde{\eta}^n_{t}=\sum_{j=1}^n\widetilde{\eta}^1_{t,j},
\end{equation}
where $\{\widetilde{\eta}^1_{t,j}\}_{1\leq j\leq n}$ are independent copies of $\widetilde{\eta}^1_{t}$. Therefore,
\[
E\left(\exp\big\{\Lambda \widetilde{\eta}^n_{t}\big\}\right)<+\infty
\]
for all $0\leq \Lambda<\log\frac{1}{1-e^{-t}}$ and all $n\geq 1$. As a result, let $\theta=\frac{1}{K_7}\log \frac{1}{1-e^{-T_0K_6K_7}}$, then Lemma \ref{lemma 3.1} follows from the above coupling relationships between $\{X_t^n\}_{t\geq 0}, \{\eta_t^n\}_{t\geq 0}$ and $\{\widetilde{\eta}_t^n\}_{t\geq 0}$.

\qed

\proof[Proof of Lemma \ref{lemma 3.2}]

By classic theory of Yule process, $E\widetilde{\eta}_t^1=e^t$. As we have introduced in the proof of Lemma \ref{lemma 3.1},
\[
E\exp\{\Lambda \widetilde{\eta}_t^1\}<+\infty
\]
for sufficiently small $\Lambda>0$. Hence, according to Equation \eqref{equ 3.1} and large deviation of the sum of i.i.d stochastic variables
(see Chapter 2 of \cite{Dembo1997}), there exists $K_8(t)>0$ such that
\[
P\left(\widetilde{\eta}_t^n> 2ne^t\right)\leq \exp\big\{-nK_8(t)\big\}
\]
for all $n\geq 1$. As we have shown in the proof of Lemma \ref{lemma 3.1},
\[
\sup_{0\leq t\leq T_0}\left\|X_t^n\right\|\leq \eta^n_{_{T_0}}=K_7\widetilde{\eta}^{\left\|X_0^n\right\|/K_7}_{_{T_0K_6K_7}}
=K_7\widetilde{\eta}^{n\left\|x_0\right\|/K_7}_{_{T_0K_6K_7}}
\]
in the sense of coupling. As a result, Lemma \ref{lemma 3.2} holds with $K_3=\frac{\left\|x_0\right\|}{K_7} K_8(T_0K_6K_7)$
and $K_2=2\left\|x_0\right\|e^{T_0K_6K_7}$.

\qed

\proof[Proof of Lemma \ref{lemma 3.3}]

By Assumptions (4) and (5), there exists $K_9\in (0,+\infty)$ such that
\[
\left\|\sum_{l\in \mathcal{A}}lF_l(x)-\sum_{l\in \mathcal{A}}lF_l(y)\right\|\leq K_9\left\|x-y\right\|
\]
for any $x,y\in G$. Then, according to Theorem 2.2 of \cite{Kurtz1978}, there exists independent Poisson processes $\{\beta_l(t):t\geq 0\}_{l\in \mathcal{A}}$ with rate one such that
\[
\sup_{0\leq t\leq T_0}\left\|\frac{X_t^n}{n}-X_t\right\|\leq A_0^n\exp\big\{K_9T_0\big\},
\]
where
\[
A_0^n=\frac{1}{n}\sup_{0\leq t\leq T_0}\Bigg\|\sum_{l\in \mathcal{A}}l\widehat{\beta}_l\left(n\int_0^tF_l\left(\frac{X_s^n}{n}\right)ds\right)\Bigg\|
\]
while $\widehat{\beta}_l(t)=\beta_l(t)-t$. Let $K_{10}=\sup_{l\in \mathcal{A}}\left\|l\right\|$. Conditioned on $\sup_{_{0\leq t\leq T_0}}\left\|X_t^n\right\|\leq nK_2$,
\[
\sup_{0\leq t\leq T_0}\Bigg\|\sum_{l\in \mathcal{A}}l\widehat{\beta}_l\left(n\int_0^tF_l\left(\frac{X_s^n}{n}\right)ds\right)\Bigg\|
\leq K_{10}\sum_{l\in \mathcal{A}}\sup_{0\leq s\leq nK_{11}T_0}\big|\widehat{\beta}_l(s)\big|,
\]
where $K_{11}=\sup_{_{l\in \mathcal{A}, \left\|x\right\|\leq K_2}}F_l(x)$. Then, according to Lemma \ref{lemma 3.2},
\begin{align}\label{equ 3.2}
&P\left(\sup_{0\leq t\leq T_0}\left\|\frac{X_t^n}{n}-X_t\right\|\geq \epsilon\right)
\leq \exp\{-K_3n\}+P\left(\sum_{l\in \mathcal{A}}\sup_{0\leq s\leq nK_{11}T_0}\big|\widehat{\beta}_l(s)\big|\geq \frac{n\epsilon e^{-K_9T_0}}{K_{10}}\right)
\notag\\
& \leq \exp\{-K_3n\}+\sum_{l\in \mathcal{A}}P\left(\sup_{0\leq s\leq nK_{11}T_0}\big|\widehat{\beta}_l(s)\big|\geq \frac{n\epsilon e^{-K_9T_0}}{K_{10}|\mathcal{A}|}\right)
\end{align}
for all $n$, where $|\mathcal{A}|$ is the cardinality of $\mathcal{A}$. By the property of Poisson process, for any $\delta>0, T_1>0$, there exists $K_{12}(\delta, T_1)\in (0,+\infty)$ such that
\begin{equation}\label{equ 3.3}
P\left(\sup_{0\leq s\leq nT_1}\big|\widehat{\beta}_l(s)\big|\geq n\delta\right)\leq \exp\big\{-nK_{12}(\delta, T_1)\big\}
\end{equation}
for sufficiently large $n$. For readers not familiar with this property, we put a proof at the end of this section. By Equations \eqref{equ 3.2} and \eqref{equ 3.3},
\[
P\left(\sup_{0\leq t\leq T_0}\left\|\frac{X_t^n}{n}-X_t\right\|\geq \epsilon\right)
\leq \exp\{-K_3n\}+|\mathcal{A}|\exp\left\{-nK_{12}\left(\frac{\epsilon e^{-K_9T_0}}{K_{10}|\mathcal{A}|},~K_{11}T_0\right)\right\}
\]
for sufficiently large $n$. As a result, Lemma \ref{lemma 3.3} holds with
\[
K_4(\epsilon)=\frac{1}{2}\min\left\{K_3, K_{12}\left(\frac{\epsilon e^{-K_9T_0}}{K_{10}|\mathcal{A}|},~K_{11}T_0\right)\right\}.
\]

\qed

At the end of this section, we give the proof of Equation \eqref{equ 3.3}.

\proof[Proof of Equation \eqref{equ 3.3}]

For simplicity, we write $\beta_l$ as $\beta$ since $\{\beta_l\}_{l\in \mathcal{A}}$ are i.i.d.. Since $\{\beta(t):t\geq 0\}$ is an independent increment process with $E\beta(t)=t$ for any $t\geq 0$, $\{\widehat{\beta}(t)=\beta(t)-t:~t\geq 0\}$ is a martingale. For any $\theta\neq 0$, $e^{\theta x}$ is a convex function with $x$, hence $\exp\{\theta\widehat{\beta}(t):t\geq 0\}$ is a submartingale. Then, by Doob's inequality,
\begin{align*}
P\left(\sup_{0\leq s\leq T_1n}\widehat{\beta}(s)\geq n\delta\right)&=P\left(\sup_{0\leq s\leq T_1n}e^{\theta\widehat{\beta}(s)}\geq e^{n\delta\theta}\right) \notag\\
&\leq e^{-n\delta\theta}Ee^{\theta\widehat{\beta}(T_1n)}= \exp\big\{-n[\delta\theta+T_1(1+\theta-e^\theta)]\big\}
\end{align*}
for any $\delta>0$ and $\theta>0$. Since $0\delta+T_1(1+0-e^0)=0$ and
\[
\frac{d}{d\theta}\left(\delta\theta+T_1(1+\theta-e^\theta)\right)\Big|_{\theta=0}=\delta>0,
\]
there exists $\theta_1>0$ such that $\delta\theta_1+T_1(1+\theta_1-e^{\theta_1})>0$ and
\[
P\left(\sup_{0\leq s\leq T_1n}\widehat{\beta}(s)\geq n\delta\right)\leq \exp\big\{-n[\delta\theta_1+T_1(1+\theta_1-e^{\theta_1})]\big\}.
\]
According to a similar analysis, there exists $\theta_2>0$ such that $ \theta_2\delta+T_1(1-\theta_2-e^{-\theta_2})>0$ and
\begin{align*}
P\left(\inf_{0\leq s\leq T_1n}\widehat{\beta}(s)\leq -n\delta\right)&=P\left(\sup_{0\leq s\leq nT_1}e^{-\theta_2\widehat{\beta}(s)}\geq e^{\theta_2n\delta}\right)\\
&\leq \exp\big\{-n[\theta_2\delta+T_1(1-\theta_2-e^{-\theta_2})]\big\}.
\end{align*}
As a result, Equation \eqref{equ 3.3} holds with
\[
K_{12}\left(\delta, T_1\right)=\frac{1}{2}\min\big\{\delta\theta_1+T_1(1+\theta_1-e^{\theta_1}), ~\theta_2\delta+T_1(1-\theta_2-e^{-\theta_2})\big\}.
\]

\qed

\section{Proof of lower bounds}\label{section four}
In this section we give the proof of the lower bound. As a preparation, we first introduce some notations and then define an exponential martingale. For each $l\in \mathcal{A}$ and $t\geq 0$, let $\xi_{t,l}^n$ be the convex combination of $X_t$ and $\frac{X_t^n}{n}$ such that
\[
F_l(\frac{X_t^n}{n})-F_l(X_t)=(\nabla F_l)(\xi_{t,l}^n)\cdot \left(\frac{X_t^n}{n}-X_t\right).
\]
Note that the existence of $\xi_{t,l}^n$ follows from Lagrange's mean value theorem. We denote by $\Omega_n$ the generator of $\{X_t^n\}_{t\geq 0}$, i.e.,
\[
\Omega_nf(x)=\sum_{l\in \mathcal{A}}nF_l\left(\frac{x}{n}\right)\big[f(x+l)-f(x)\big]
\]
for any sufficiently smooth $f: \mathbb{R}^d\rightarrow \mathbb{R}$. For any $f_1, f_2\in C^{2,1}\left([0,T_0)\times \mathbb{R}^d\right)$,
let
\[
\mathcal{M}_{f_1}^n(t)=f_1(t,X_t^n)-f_1(0,X_0^n)-\int_0^t\left(\frac{\partial}{\partial s}+\Omega_n\right)f(s,X_s^n) ds,
\]
and
\begin{align*}
\mathcal{N}^n_{f_1,f_2}(t)&=\mathcal{M}_{f_1}^n(t)\mathcal{M}_{f_2}^n(t)\\
&-\int_0^t\Omega_n(f_1f_2)(s,X_s^n)-f_1(s,X_s^n)\Omega_nf_2(s,X_s^n)-f_2(s,X_s^n)\Omega_nf_1(s,X_s^n)ds,
\end{align*}
then according to properties of continuous-time Markov processes (see Section 5 of Appendix 1 of \cite{Kipnis1999}), $\{\mathcal{M}^n_{f_1}(t)\}_{0\leq t \leq T_0}$ and
$\{\mathcal{N}^n_{f_1,f_2}(t)\}_{0\leq t\leq T_0}$ are both martingales. That is to say,
\begin{equation}\label{equ 4.cross variation}
d\langle\mathcal{M}^n_{f_1}, \mathcal{M}^n_{f_2}\rangle_t=\left(\Omega_n(f_1f_2)-f_1\Omega_nf_2-f_2\Omega_nf_1\right)dt.
\end{equation}
Note that in this paper $\langle\cdot\rangle$ and $[\cdot]$ are defined in the same way as that defined in \cite{Schuppen1974}, i.e, for a local martingale $M$, $\langle M\rangle$ is the unique predictable increasing process such that $M^2-\langle M\rangle$ is a local martingale while $[M]$ is the quadratic-variation process of $M$ (which is not equal to $\langle M\rangle$ when $M$ is not continuous). For two local martingales $M_1, M_2$, $\langle M_1, M_2\rangle$ and $[M_1, M_2]$ are defined as
\[
\langle M_1, M_2\rangle=\frac{\langle M_1+M_2\rangle-\langle M_1-M_2\rangle}{4} \text{~and~} [M_1, M_2]=\frac{[M_1+M_2]-[M_1-M_2]}{4}.
\]

To utilize above martingales, for any $g\in C^2\left([0,T_0], \mathbb{R}^d\right)$, let
$f_{n,g}(t,x)=\frac{a_n}{n}g_t\cdot (x-nX_t)$ and consequently
\[
f_{n,g}(t,X_t^n)=\frac{a_n}{n}g_t\cdot(X_t^n-nX_t),
\]
then by direct calculation and Equation \eqref{equ 4.cross variation},
\begin{align}\label{equ 4. Ztg}
df_{n,g}(t,X_t^n)&=\left(\frac{\partial}{\partial t}+\Omega_n\right)f_{n,g}(t, X_t^n)dt+dM_t(f_{n,g})\\
&=\frac{a_n}{n}g_t^\prime\cdot(X_t^n-nX_t)dt+\frac{a_n}{n}g_t^{\mathsf{T}}\sum_{l\in \mathcal{A}}l(\nabla^{\mathsf{T}}F_l)(\xi_{t,l}^n)(X_t^n-nX_t)dt+dM_t(f_{n,g}), \notag
\end{align}
where $\{M_t(f_{n,g})\}_{0\leq t\leq T_0}$ is a martingale with $M_0(f_{n,g})=0$ and
\[
d\left\langle M(f_{n,g}), M(f_{n,h})\right\rangle_t=\frac{a_n^2}{n}g_t^{\mathsf{T}}\big[\sum_{l\in \mathcal{A}}lF_l(\frac{X_t^n}{n})l^{\mathsf{T}}\big]h_tdt
\]
for any $g,h\in C^2\left([0,T_0], \mathbb{R}^d\right)$. For later use, we define
\[
H_{n,g}(t,x)=\exp\{f_{n,g}(t,x)\}
\]
and consequently $H_{n,g}(t, X_t^n)=\exp\{\frac{a_n}{n}g_t\cdot(X_t^n-nX_t)\}$ for any $t\geq 0$. Our exponential martingale is defined according to the following lemma.
\begin{lemma}\label{lemma 4.1 exponential martingale}
For any $g\in C^2\left([0,T_0], \mathbb{R}^d\right)$, let
\[
\omega_t^n(g)=\frac{H_{n,g}(t, X_t^n)}{H_{n,g}(0, X_0^n)}\exp\left\{-\int_0^t \frac{(\frac{\partial}{\partial s}+\Omega_n)H_{n,g}(s,X_s^n)}{H_{n,g}(s, X_s^n)}ds\right\},
\]
then there exists $N(g)\geq 1$ such that
$\{\omega_t^n(g)\}_{0\leq t\leq T_0}$ is a martingale with expectation $1$ for each $n\geq N(g)$.
\end{lemma}

\proof

According to Integration-by-parts formula (see Volume 2, Chapter 6, section 38 of \cite{Rogers1986}) and direct calculation,
\begin{equation}\label{equ revision 4.1}
d\omega_t^n(g)=\Lambda^n_{t-}(g)dM_t(H_{n,g}),
\end{equation}
where
\[
M_t(H_{n,g})=H_{n,g}(t, X_t^n)-H_{n,g}(0, X_0^n)-\int_0^t \left(\frac{\partial}{\partial s}+\Omega_n\right)H_{n,g}(s, X_s^n)ds
\]
and
\[
\Lambda^n_t(g)=\frac{1}{H_{n,g}(0, X_0^n)}\exp\left\{-\int_0^t \frac{(\frac{\partial}{\partial s}+\Omega_n)H_{n,g}(s, X_s^n)}{H_{n,g}(s, X_s^n)}ds\right\}.
\]
As we have recalled, $\{M_t(H_{n,g})\}_{0\leq t\leq T_0}$ is a martingale, hence $\{\omega_t^n(g)\}_{0\leq t\leq T_0}$ is a local martingale. Therefore, to check that $\{\omega_t^n(g)\}_{0\leq t\leq T_0}$ is a martingale for large $n$, we only need to show that
\begin{equation}\label{equ revision 4.2}
E\Bigg(\int_0^{T_0}\left(\Lambda^n_{s-}(g)\right)^2d[M(H_{n,g})]_s\Bigg)<+\infty
\end{equation}
for sufficiently large $n$. By direct calculation and Taylor's expansion formula up to the second order,
\begin{align*}
\Lambda^n_t(g)=&\exp\Bigg\{-\frac{a_n}{n}\int_0^tg_s^\prime\cdot (X_s^n-nX_s)ds+a_n\int_0^tg_s\cdot \left(\sum_llF_l(X_s)\right)ds\\
&\quad\quad\quad\quad-\int_0^t\sum_lnF_l(\frac{X_s^n}{n})\left(e^{\frac{a_n}{n}g_s\cdot l}-1\right)ds\Bigg\}\\
=&\exp\Bigg\{-\frac{a_n}{n}\int_0^tg_s^\prime\cdot (X_s^n-nX_s)ds-\frac{a_n}{n}\int_0^tg_s\cdot\left(\sum_ll(\nabla^{\mathsf{T}}F_l)(\xi_{s,l}^n) (X_s^n-nX_s)\right)ds\\
&\quad\quad\quad\quad-\frac{a_n^2}{2n}(1+o(1))\int_0^tg_s^{\mathsf{T}}\left(\sum_llF_l(\frac{X_s^n}{n})l^{\mathsf{T}}\right)g_sds\Bigg\}.
\end{align*}
According to the definition of the quadratic-variation process of a discontinuous martingale (see Section 2 of \cite{Schuppen1974}),
\[
[M(H_{n,g})]_t=\sum_{s\leq t}\left(H_{n,g}(s, X_s^n)-H_{n,g}(s-, X_{s-}^n)\right)^2.
\]
As a result, according to Assumptions (3)-(5) and the coupling relationship given in the proof of Lemma \ref{lemma 3.1}, there exists $K_{13}\in (0,+\infty)$ depending on $T_0, g$ and $\mathcal{A}$ such that
\[
\sup_{0\leq t\leq T_0}\left(\Lambda^n_t(g)\right)^2\leq \exp\left\{a_nK_{13}+\frac{a_n}{n}K_{13}\sup_{0\leq t\leq T_0}\left\|X_t^n\right\|+\frac{a_n^2}{n^2}K_{13}\sup_{0\leq t\leq T_0}\left\|X_t^n\right\|\right\}
\]
and
\[
\sup_{0\leq t\leq T_0}[M(H_{n,g})]_t\leq \widetilde{\eta}_{_{K_{13}T_0}}^{n K_{13}\left\|x_0\right\|}\exp\left\{a_nK_{13}+K_{13}\frac{a_n}{n}\widetilde{\eta}_{_{K_{13}T_0}}^{n K_{13}\left\|x_0\right\|}\right\}
\]
for sufficiently large $n$, where $\{\widetilde{\eta}_t^n\}_{t\geq 0}$ is the Yule process defined as in Section \ref{section three}. Then, Equation \eqref{equ revision 4.2} follows directly from Lemma \ref{lemma 3.1} and the facts that $\frac{a_n}{n}\rightarrow 0$ while $\widetilde{\eta}_t^1$ follows a Geometric distribution.

\qed

Let $P$ be the probability measure of our DDMC, then for $g\in C^2\left([0,T_0], \mathbb{R}^d\right)$ and each $n\geq N(g)$, let $P_n^g$ be the probability measure such that
\[
\frac{dP_n^g}{dP}=\omega^n_{_{T_0}}(g),
\]
then we have the following laws of large numbers.

\begin{lemma}\label{lemma 4.2 LLNunderPn}
As $n\rightarrow+\infty$, $\{\frac{X_t^n}{n}\}_{0\leq t\leq T_0}$ converges in $P_n^g$-probability  to $\{X_t\}_{0\leq t\leq T_0}$, where $\{X_t\}_{0\leq t\leq T_0}$ is defined as in Equation \eqref{equ 1.1 ODE}.
\end{lemma}

\begin{lemma}\label{lemma 4.3 LLNunderGirsanov}
As $n\rightarrow+\infty$, $\big\{\frac{X_t^n-nX_t}{a_n}\big\}_{0\leq t\leq T_0}$ converges in $P_n^g$-probability  to the solution of the ODE
\begin{equation}\label{equ 4.1.5 ODEGirsanov}
\begin{cases}
&\frac{d}{dt}y_t=b_ty_t+\sigma_tg_t \text{~for~}0\leq t\leq T_0, \\
&y_0=0,
\end{cases}
\end{equation}
where $b_t$ and $\sigma_t$ are defined as in Section \ref{section one}.
\end{lemma}

\proof[Proof of Lemma \ref{lemma 4.2 LLNunderPn}]

For any $\epsilon>0$, according to Cauchy-Schwartz's inequality,
\begin{align}\label{equ 4.1}
P_n^g\left(\sup_{0\leq t\leq T_0}\left\|\frac{X_t^n}{n}-X_t\right\|\geq \epsilon\right)&=E\left(\omega^n_{_{T_0}}(g)1_{_{\left\{\sup_{0\leq t\leq T_0}\left\|\frac{X_t^n}{n}-X_t\right\|\geq \epsilon\right\}}}\right)\\
&\leq \sqrt{E\left(\left(\omega^n_{_{T_0}}(g)\right)^2\right)}\sqrt{P\left(\sup_{0\leq t\leq T_0}\left\|\frac{X_t^n}{n}-X_t\right\|\geq \epsilon\right)}.\notag
\end{align}
According to the definitions of $\omega^n_t(g), M_t(f_{n,g})$ and direct calculation,
\begin{align}\label{equ revision 4.calculationofomega}
\omega^n_{t}(g)&=\Lambda^n_{t}(g)\exp\left\{\frac{a_n}{n}g_{_{t}}\cdot(X_{t}^n-nX_{t})-\frac{a_n}{n}g_0\cdot(X_0^n-nX_0)\right\}\notag\\
&=\exp\Bigg\{\frac{a_n}{n}g_{_{t}}\cdot(X_{t}^n-nX_{t})-\frac{a_n}{n}g_0\cdot(X_0^n-nX_0)\notag\\
&\quad\quad-\frac{a_n}{n}\int_0^{t}g_s^\prime\cdot (X_s^n-nX_s)ds-\frac{a_n}{n}\int_0^{t}g_s\cdot\left(\sum_ll(\nabla^{\mathsf{T}}F_l)(\xi_{s,l}^n) (X_s^n-nX_s)\right)ds\notag\\
&\quad\quad-\frac{a_n^2}{2n}(1+o(1))\int_0^{t}g_s^{\mathsf{T}}\left(\sum_llF_l(\frac{X_s^n}{n})l^{\mathsf{T}}\right)g_sds\Bigg\}\notag\\
&=\exp\left\{M_{t}(f_{n,g})-\frac{1}{2}(1+o(1))\left\langle M(f_{n,g})\right\rangle_{t}\right\}.
\end{align}

Then, by Assumptions (4) and (5), there exists $K_{14}\in (0,+\infty)$ depending on $g$ and $T_0$ such that
\[
\omega^n_{_{T_0}}(g)\leq \exp\Bigg\{K_{14}\left(a_n\sup_{0\leq t\leq T_0}\left\|\frac{X_t^n}{n}-X_t\right\|+\frac{a_n^2}{n^2}\sup_{0\leq t\leq T_0}\left\|X_t^n\right\|\right)\Bigg\}.
\]
Let $K_{15}=\sup_{0\leq t\leq T_0}\left\|X_t\right\|$, since $\frac{a_n}{n}\rightarrow 0$,
\begin{align*}
\omega^n_{_{T_0}}(g)&\leq \exp\left\{K_{14}\frac{a_n}{n}\sup_{0\leq t\leq T_0}\left\|X_t^n\right\|+a_nK_{14}K_{15}+K_{14}\frac{a_n^2}{n^2}\sup_{0\leq t\leq T_0}\left\|X_t^n\right\|\right\}\\
&\leq \exp\left\{\frac{2K_{14}a_n}{n}\sup_{0\leq t\leq T_0}\left\|X_t^n\right\|+a_nK_{14}K_{15}\right\}\\
&=\exp\left\{\frac{K_{16}a_n}{n}\sup_{0\leq t\leq T_0}\left\|X_t^n\right\|+a_nK_{17}\right\}
\end{align*}
for sufficiently large $n$, where $K_{16}=2K_{14}$ and $K_{17}=K_{14}K_{15}$. Therefore,
\[
E\left(\left(\omega^n_{_{T_0}}(g)\right)^2\right)\leq e^{2a_nK_{17}}E\exp\left\{\frac{2K_{16}a_n}{n}\sup_{0\leq t\leq T_0}\left\|X_t^n\right\|\right\}.
\]
As we have introduced in the proof of Lemma \ref{lemma 3.2},
\[
\sup_{0\leq t\leq T_0}\left\|X_t^n\right\|\leq K_7\widetilde{\eta}^{n\left\|x_0\right\|/K_7}_{_{T_0K_6K_7}}
\]
while $\widetilde{\eta}^1_t$ follows geometric distribution with rate $e^{-t}$ and $\widetilde{\eta}^n_t$ is the sum of $n$ i.i.d. copies of $\widetilde{\eta}^1_t$. Therefore, according to the fact that $\frac{a_n}{n}\rightarrow 0$,
\begin{align*}
E\exp\left\{\frac{2K_{16}a_n}{n}\sup_{0\leq t\leq T_0}\left\|X_t^n\right\|\right\}& \leq \Bigg(\frac{e^{\frac{2K_{16}K_7a_n}{n}}e^{-T_0K_6K_7}}{1-(1-e^{-T_0K_6K_7})e^{\frac{2K_{16}K_7a_n}{n}}}\Bigg)^{n\left\|x_0\right\|/K_7}\\
&=e^{a_n\left(K_{18}+o(1)\right)}
\end{align*}
for sufficiently large $n$, where $K_{18}=2\frac{K_{16}}{K_7}\left\|x_0\right\| e^{T_0K_6K_7}$. As a result, for sufficiently large $n$,
\[
E\left(\left(\omega^n_{_{T_0}}(g)\right)^2\right)\leq e^{a_n(K_{18}+2K_{17}+o(1))}
\]
and hence
\[
P_n^g\left(\sup_{0\leq t\leq T_0}\left\|\frac{X_t^n}{n}-X_t\right\|\geq \epsilon\right)\leq e^{a_n(\frac{K_{18}}{2}+K_{17}+o(1))}e^{-\frac{K_4(\epsilon)n}{2}}
\]
by Lemma \ref{lemma 3.3} and Equation \eqref{equ 4.1}. Since $\frac{a_n}{n}\rightarrow0$,
\[
\lim_{n\rightarrow+\infty}P_n^g\left(\sup_{0\leq t\leq T_0}\left\|\frac{X_t^n}{n}-X_t\right\|\geq \epsilon\right)=0
\]
for any $\epsilon>0$ and hence Lemma \ref{lemma 4.2 LLNunderPn} holds.

\qed

\proof[Proof of Lemma \ref{lemma 4.3 LLNunderGirsanov}]

For $1\leq i\leq d$, let $e_i$ be the $i$th elementary unit vector of $\mathbb{R}^d$, i.e.,
\[
e_i=(0,\ldots,0,\mathop 1\limits_{i \text{th}},0,\ldots,0)^{\mathsf{T}}
\]
and denote $\frac{a_n}{n}e_i\cdot (x-nX_t)$ by $f_{n,i}(t,x)$,
then by Equation \eqref{equ 4. Ztg}, $f_{n,i}(t, X_t^n)=\frac{a_n}{n}e_i\cdot (X_t^n-nX_t)$ satisfies
\[
df_{n,i}(t, X_t^n)=\frac{a_n}{n}e_i^{\mathsf{T}}\sum_{l\in \mathcal{A}}l(\nabla^{\mathsf{T}}F_l)(\xi_{t,l}^n)(X_t^n-nX_t)dt+dM_t(f_{n,i}),
\]
where $\{M_t(f_{n,i})\}_{0\leq t\leq T_0}$ is a martingale for each $1\leq i\leq d$ and
\[
d\langle M(f_{n,i}), M(f_{n,j})\rangle_t=\frac{a_n^2}{n}\left[\sum_{l\in \mathcal{A}}lF_l\left(\frac{X_t^n}{n}\right)l^{\mathsf{T}}\right]_{ij}dt.
\]
Since $\frac{a_n}{n}(X_t^n-nX_t)=\left(f_{n,1}(t, X_t^n),\ldots, f_{n,d}(t, X_t^n)\right)^{\mathsf{T}}$,
\[
d\frac{a_n}{n}(X_t^n-nX_t)=\frac{a_n}{n}\sum_{l\in \mathcal{A}}l(\nabla^{\mathsf{T}}F_l)(\xi_{t,l}^n)(X_t^n-nX_t)dt+d\mathcal{M}^n_t,
\]
where $\mathcal{M}^n_t=\left(M_t(f_{n,1}), \ldots, M_t(f_{n,d})\right)^{\mathsf{T}}$. Then,
\begin{equation}\label{equ revision 4.3 two}
d\frac{X_t-nX_t}{a_n}=\sum_{l\in \mathcal{A}}l(\nabla^{\mathsf{T}}F_l)(\xi_{t,l}^n)\frac{X_t^n-nX_t}{a_n}dt+\frac{n}{a_n^2}d\mathcal{M}^n_t.
\end{equation}
Let $M_t(H_{n,g})$ be defined as in the proof of Lemma \ref{lemma 4.1 exponential martingale}, then we define
\[
\widetilde{M}_t(H_{n,g})=\int_0^t \frac{1}{H_{n,g}(s-, X_{s-}^n)}dM_s(H_{n,g})
\]
for each $0\leq t\leq T_0$. Then, by Equation \eqref{equ revision 4.1},
\begin{equation}\label{equ revision 4.3}
d\omega_t^n(g)=\Lambda^n_{t-}dM_t(H_{n,g})=\Lambda^n_{t-}H_{n,g}(t-, X_{t-}^n)d\widetilde{M}_t(H_{n,g})=\omega_{t-}^n(g)d\widetilde{M}_t(H_{n,g}).
\end{equation}
For each $0\leq t\leq T_0$ and $1\leq i\leq d$, let
\[
\widehat{M}_t(f_{n,i})=M_t(f_{n,i})-\langle M(f_{n,i}), \widetilde{M}(H_{n,g})\rangle_t,
\]
then by Equation \eqref{equ revision 4.3} and Theorem 3.2 of \cite{Schuppen1974}, which is a generalized version of Girsanov's Theorem, $\{\widehat{M}_t(f_{n,i})\}_{0\leq t\leq T_0}$ is a local martingale under $P_n^g$ and
\[
[\widehat{M}(f_{n,i})]_t=[M(f_{n,i})]_t
\]
under $P$ and $P_n^g$. Note that by Equation \eqref{equ 4.cross variation} and direct calculation,
\begin{align}\label{equ revision 4.4}
d\langle M(f_{n,i}), \widetilde{M}(H_{n,g})\rangle_t&=\frac{1}{H_{n,g}(t-, X_{t-}^n)}d\langle M(f_{n,i}), M(H_{n,g})\rangle_t \\
&=\frac{a_n}{n}\sum_{l}nF_l(\frac{X_t^n}{n})e_i\cdot l\left(e^{\frac{a_n}{n}g_t\cdot l}-1\right)dt \notag\\ &=\frac{a_n^2}{n}e_i^{\mathsf{T}}\Big[\sum_llF_l(\frac{X_t^n}{n})l^{\mathsf{T}}g_t(1+o(1))\Big]dt. \notag
\end{align}
Let $\widehat{\mathcal{M}}^n_t=\left(\widehat{M}_t(f_{n,1}),\ldots,\widehat{M}_t(f_{n,d})\right)^{\mathsf{T}}$, then by Equations \eqref{equ revision 4.3 two} and \eqref{equ revision 4.4},
\begin{equation}\label{equ revision 4.5}
d\frac{X_t-nX_t}{a_n}=b_t^n\frac{X_t^n-nX_t}{a_n}dt+\sigma_t^ng_tdt+\frac{n}{a_n^2}d\widehat{\mathcal{M}}_t^n,
\end{equation}
where $b_t^n=\sum_{l\in \mathcal{A}}l(\nabla^{\mathsf{T}}F_l)(\xi_{t,l}^n)$ and $\sigma_t^n=\sum_llF_l(\frac{X_t^n}{n})
l^{\mathsf{T}}(1+o(1))$.

By Lemma \ref{lemma 4.2 LLNunderPn}, $\left\{\sum_{l\in \mathcal{A}}l(\nabla^{\mathsf{T}}F_l)(\xi_{t,l}^n)\right\}_{0\leq t\leq T_0}$ converges in $P_n^g$-probability to
\[
\left\{\sum_{l\in \mathcal{A}}l(\nabla^{\mathsf{T}}F_l)(X_t)\right\}_{0\leq t\leq T_0}=\{b_t\}_{0\leq t\leq T_0}
\]
and $\left\{\sum_{l\in \mathcal{A}}lF_l(\frac{X_t^n}{n})l^{\mathsf{T}}\right\}_{0\leq t\leq T_0}$ converges in $P_n^g$-probability to
\[
\left\{\sum_{l\in \mathcal{A}}lF_l(X_t)l^{\mathsf{T}}\right\}_{0\leq t\leq T_0}=\{\sigma_t\}_{0\leq t\leq T_0}
\]
as $n\rightarrow+\infty$.

By Assumption (5), there exists $K_{20}\in (0, +\infty)$ such that
\[
\sup_{0\leq t\leq T_0}\left\|\sum_{l\in \mathcal{A}}l(\nabla^{\mathsf{T}}F_l)(\xi_{t,l}^n)x\right\|\leq K_{20}\left\|x\right\|
\]
for any $x\in \mathbb{R}^d$. Consequently, by Grownwall's inequality,
\[
\left\|\frac{X_t-nX_t}{a_n}-y_t\right\|\leq \varpi_n e^{K_{20}T_0}
\]
for any $0\leq t\leq T_0$, where
\[
\varpi_n=\sup_{0\leq t\leq T_0}\left(\int_0^{T_0}\left\|(b_s^n-b_s)y_s\right\|ds+\int_0^{T_0} \left\|(\sigma_s^n-\sigma_s)g_s\right\|ds+\frac{n}{a_n^2}\left\|\widehat{\mathcal{M}}_t^n\right\|\right).
\]

As we have shown, $b_s^n-b_s, \sigma_s^n-\sigma_s$ converges in $P_n^g$-probability to $0$. Hence, to complete this proof, we only need to show that $\sup_{0\leq t\leq T}\frac{n}{a_n^2}\left\|\widehat{\mathcal{M}}_t^n\right\|$ converges in $P_n^g$-probability to $0$, i.e., $\sup\limits_{0\leq t\leq T_0}\frac{n}{a_n^2}|\widehat{M}_t(f_{n,i})|$ converges in $P_n^g$-probability to $0$ for each $1\leq i\leq d$.

As we have recalled, under $P_n^g$,
\begin{equation}\label{equ revision 4.6}
[\widehat{M}(f_{n,i})]_t=[M(f_{n,i})]_t=\sum_{s\leq t}\left(\frac{a_ne_i\cdot \left(X_s^n-X_{s-}^n\right)}{n}\right)^2
\end{equation}
according to the generalized Girsanov's theorem introduced in \cite{Schuppen1974}. For any $\delta>0$, let $\tau_n(\delta)=\inf\{t:~[\widehat{M}(f_{n,i})]_t\geq \delta\}$, then by Equation \eqref{equ revision 4.6} and Assumption (4), there exists $K_{21}\in (0,+\infty)$ depending on $i$ such that
\[
[\widehat{M}(f_{n,i})]_{\tau_n(\delta)}\leq \delta+\frac{a_n^2K_{21}}{n^2}
\]
for any $\delta>0$. Then, by Doob's inequality,
\begin{align*}
&P_n^g\left(\sup_{0\leq t\leq T_0}\frac{n}{a_n^2}|\widehat{M}_t(f_{n,i})| \geq \epsilon\right)=P_n^g\left(\sup_{0\leq t\leq T_0}\frac{n}{a_n^2}|\widehat{M}_t(f_{n,i})| \geq \epsilon, \tau_n\left(\frac{a_n^4}{n^2}\delta\right)>T_0\right)\\
&\quad\quad+P_n^g\left(\sup_{0\leq t\leq T_0}\frac{n}{a_n^2}|\widehat{M}_t(f_{n,i})| \geq \epsilon, \tau_n\left(\frac{a_n^4}{n^2}\delta\right)\leq T_0\right)\\
&\leq \frac{\frac{a_n^4}{n^2}\delta+\frac{a_n^2}{n^2}K_{21}}{\frac{a_n^4}{n^2}\epsilon^2}+P_n^g\left(\tau_n\left(\frac{a_n^4}{n^2}\delta\right)\leq T_0\right)=\frac{\delta+\frac{1}{a_n^2}K_{21}}{\epsilon^2}+P_n^g\left([\widehat{M}(f_{n,i})]_{T_0}\geq \frac{a_n^4}{n^2}\delta\right)
\end{align*}
for any $\epsilon, \delta>0$. Consequently, we only need to show that
\begin{equation}\label{equ revision 4.7}
\lim_{n\rightarrow+\infty}P_n^g\left([\widehat{M}(f_{n,i})]_{T_0}\geq \frac{a_n^4}{n^2}\delta\right)=0
\end{equation}
for any $\delta>0$ to finish this proof. To prove Equation \eqref{equ revision 4.7}, we let
\[
\Xi^n=\sum_{t\leq T_0}1_{\{X^n_t\neq X^n_{t-}\}},
\]
i.e., $\Xi^n$ is the number of jumps in $\{X_t^n\}_{0\leq t\leq T_0}$. Then by Equation \eqref{equ revision 4.6} and Assumption (4), there exists $K_{22}\in (0,+\infty)$ depending on $i$ such that $[\widehat{M}(f_{n,i})]_{T_0}\leq \frac{a_n^2K_{22}\Xi^n}{n^2}$. Therefore, according to a similar analysis with that
in the proof of Lemma \ref{lemma 4.2 LLNunderPn} and Cauchy-Schwartz's inequality,
\begin{align}\label{equ revision 4.8}
P_n^g\left([\widehat{M}(f_{n,i})]_{T_0}\geq \frac{a_n^4}{n^2}\delta\right)
\leq e^{a_n(\frac{K_{18}}{2}+K_{17}+o(1))}\sqrt{P\left(\Xi^n\geq \frac{a_n^2\delta}{K_{22}}\right)}.
\end{align}
According to a similar analysis with that in the proof of Lemma \ref{lemma 3.1}, there exists $K_{23}\in (0,+\infty)$ such that $\Xi^n$ is stochastic dominated from above by $K_{23}\widetilde{\eta}^{n K_{23}\left\|x_0\right\|}_{K_{23}T_0}$ under $P$, where $\{\widetilde{\eta}^n_t\}_{t\geq 0}$ is the Yule process defined as in Section \ref{section three}. As we have recalled, $\widetilde{\eta}^n_t$ is the sum of $n$ i.i.d copies of a random variable following a Geometric distribution. Therefore, according to the large deviation principle for the sum of i.i.d. random variables and the fact that $\frac{a_n^2}{n}\rightarrow+\infty$, there exists $K_{24}\in (0,+\infty)$ such that
\begin{equation}\label{equ revision 4.9}
P\left(\Xi^n\geq \frac{a_n^2\delta}{K_{22}}\right)\leq P\left(\frac{\widetilde{\eta}_{K_{23}T_0}^{n K_{23}\left\|x_0\right\|}}{n}\geq \frac{a_n^2\delta}{nK_{22}K_{23}}\right)\leq e^{-K_{24}n}
\end{equation}
for sufficiently large $n$. Since $\frac{a_n}{n}\rightarrow 0$, Equation \eqref{equ revision 4.7} follows from Equations \eqref{equ revision 4.8} and \eqref{equ revision 4.9} directly and the proof is complete.

\qed

To give the proof of the lower bound, we need the following lemma, which is a generalized version of Equation \eqref{equ 2.2} under the case where $\sigma_t$ is not invertible.

\begin{lemma}\label{lemma 4.4}
If $f\in \mathcal{D}([0,T_0], \mathbb{R}^d)$ makes $I(f)<+\infty$, then $f$ is absolutely continuous and there exists $\psi: [0, T_0]\rightarrow \mathbb{R}^d$ such that
\[
f^\prime_t=b_tf_t+\sigma_t\psi_t
\]
for $0\leq t\leq T_0$ and $I(f)=\frac{1}{2}\int_0^{T_0}\psi_s^{\mathsf{T}}\sigma_s\psi_s ds$.
\end{lemma}

The proof of Lemma \ref{lemma 4.4} follows from a similar strategy with that of Equation \eqref{equ 2.2}.

\proof[Proof of Lemma \ref{lemma 4.4}] For $f$ making $I(f)<+\infty$ and $g\in C^2\left([0,T_0],\mathbb{R}^d\right)$, if $\mathcal{L}_2(g)=0$ but $\mathcal{L}_{1,f}(g)\neq 0$, then
\[
I(f)\geq \sup\{\mathcal{L}_{1,f}(cg)-\frac{1}{2}\mathcal{L}_2(cg):c\in \mathbb{R}\}=\sup\{c\mathcal{L}_{1,f}(g):c\in \mathbb{R}\}=+\infty,
\]
which is contradictory. Hence, $\mathcal{L}_2(g)=0$ implies that $\mathcal{L}_{1,f}(g)=0$. For $g$ making $\mathcal{L}_2(g)\neq 0$,
$\mathcal{L}_{1,f}(cg)-\frac{1}{2}\mathcal{L}_2(cg)$ get maximum $\frac{\left(\mathcal{L}_{1,f}(g)\right)^2}{2\mathcal{L}_2(g)}$ when
$c=\frac{\mathcal{L}_{1,f}(g)}{\mathcal{L}_2(g)}$, hence
\[
I(f)=\sup\left\{\frac{\left(\mathcal{L}_{1,f}(g)\right)^2}{2\mathcal{L}_2(g)}:~g\text{~makes~}\mathcal{L}_2(g)\neq 0\right\}.
\]
Note that $\mathcal{L}_2(g)=0$ when and only when $\sigma_t^{\frac{1}{2}}g_t=0$ almost everywhere for $0\leq t\leq T_0$ (square root of  $\sigma_t$ can be defined since $\sigma_t$ is positive semi-definite). For $g,h\in L^2_\sigma\left([0,T_0], \mathbb{R}^d\right)$, we write $g\simeq h$ when and only when
\[
\sigma_t^{\frac{1}{2}}(g_t-h_t)=0 \text{~a.e..}
\]
Then, $\simeq$ is an equivalence relation. We define $[g]=\{h:h\simeq g\}$ and
\[
L^2_{\sigma, \simeq}\left([0,T_0], \mathbb{R}^d\right)=\big\{[g]:~g\in L^2_\sigma\left([0,T_0], \mathbb{R}^d\right)\big\}.
\]
For $[g], [h]\in L^2_{\sigma, \simeq}\left([0,T_0], \mathbb{R}^d\right)$, we define
\[
\langle[g], [h]\rangle_{L^2_{\sigma, \simeq}}=\int_0^{T_0} g_s^{\mathsf{T}}\sigma_sh_sds.
\]
It is easy to check that $\langle \cdot,\cdot\rangle_{L^2_{\sigma, \simeq}}$ is well-defined and $L^2_{\sigma, \simeq}\left([0,T_0], \mathbb{R}^d\right)$ is a Hilbert space under $\langle \cdot,\cdot\rangle_{L^2_{\sigma, \simeq}}$. We define
$C^2_{\sigma, \simeq}\left([0,T_0], \mathbb{R}^d\right)=\big\{[g]:~g\in C^2\left([0,T_0], \mathbb{R}^d\right)\big\}$. For $g\in C^2\left([0,T_0], \mathbb{R}^d\right)$, we define
\[
\mathcal{L}_{1,f,\simeq}([g])=\mathcal{L}_{1,f}(g).
\]
According to the fact that $g\simeq h$ implies $\mathcal{L}_2(g-h)=0$ and hence $\mathcal{L}_{1,f}(g-h)=0$, $\mathcal{L}_{1,f,\simeq}$ is well-defined and
is a linear operator on $C^2_{\sigma, \simeq}\left([0,T_0], \mathbb{R}^d\right)$. Then,
\[
I(f)=\sup\left\{\frac{\left(\mathcal{L}_{1,f, \simeq}([g])\right)^2}{2\langle[g], [g]\rangle_{L^2_{\sigma, \simeq}}}:~[g]\neq [0], [g]\in C^2_{\sigma, \simeq}\left([0,T_0], \mathbb{R}^d\right) \right\}.
\]
Since $I(f)<+\infty$ and $C^2_{\sigma, \simeq}\left([0,T_0], \mathbb{R}^d\right)$ is dense in $L^2_{\sigma, \simeq}\left([0,T_0], \mathbb{R}^d\right)$, $\mathcal{L}_{1,f,\simeq}$ can be extended to a bounded linear operator on $L^2_{\sigma, \simeq}\left([0,T_0], \mathbb{R}^d\right)$ and hence there exists
$[\psi]\in L^2_{\sigma, \simeq}\left([0,T_0], \mathbb{R}^d\right)$ such that
\[
\mathcal{L}_{1,f,\simeq}([g])=\langle [g], [\psi] \rangle_{L^2_{\sigma,\simeq}}
\]
for any $[g]\in C^2_{\sigma, \simeq}\left([0,T_0], \mathbb{R}^d\right)$ according to Riesz representation theorem. As a result,
\[
f(T_0)\cdot g(T_0)-f(0)\cdot g(0)-\int_0^{T_0}f_s\cdot g^\prime_sds=\int_0^{T_0}(b_sf_s+\sigma_s\psi_s)\cdot g_sds
\]
for any $g\in C^2\left([0,T_0], \mathbb{R}^d\right)$. Therefore, $f$ is absolutely continuous and
\[
f_t^\prime=b_tf_t+\sigma_t\psi_t.
\]
$I(f)\leq \frac{1}{2}\langle [\psi], [\psi]\rangle_{L^2_{\sigma, \simeq }}=\frac{1}{2}\int_0^{T_0}\psi_s^{\mathsf{T}}\sigma_s\psi_sds$ follows from the fact that
\[
\left(\mathcal{L}_{1,f, \simeq}([g])\right)^2=\left(\langle [g], [\psi] \rangle_{L^2_{\sigma,\simeq}}\right)^2\leq
\langle [g], [g] \rangle_{L^2_{\sigma,\simeq}}\langle [\psi], [\psi] \rangle_{L^2_{\sigma,\simeq}}
\]
for any $g\in C^2\left([0,T_0], \mathbb{R}^d\right)$ by Cauchy Schwartz's inequality. To prove $I(f)\geq \frac{1}{2}\langle [\psi], [\psi]\rangle_{L^2_{\sigma, \simeq }}$, we choose $g_n\in C^2\left([0,T_0], \mathbb{R}^d\right)$ such that $[g_n]$ converges to $[\psi]$ under the distance generated by $\langle \cdot,\cdot\rangle_{L^2_{\sigma,\simeq}}$ as $n$ grows to infinity. Then,
\[
I(f)\geq \lim_{n\rightarrow+\infty} \frac{\left(\langle[g_n], [\psi]\rangle_{L^2_{\sigma, \simeq}}\right)^2}{2\langle[g_n], [g_n]\rangle_{L^2_{\sigma, \simeq}}}
=\frac{1}{2}\langle [\psi], [\psi]\rangle_{L^2_{\sigma, \simeq }}.
\]

\qed

At the end of this section, we give the proof of the lower bound.

\proof[Proof of the lower bound]

For given open set $O\subseteq \mathcal{D}\left([0,T_0], \mathbb{R}^d\right)$, if $\inf_{f\in O}I(f)=+\infty$, then
\[
\liminf_{n\rightarrow+\infty}\frac{n}{a_n^2}\log P\left(\vartheta^n\in O\right)\geq -\inf_{f\in O}I(f)
\]
holds trivial. If $\inf_{f\in O}I(f)<+\infty$, then for any $\epsilon>0$, there exists $f_\epsilon\in O$ such that $I(f_\epsilon)\leq \inf_{f\in O}I(f)+\epsilon$. By Lemma \ref{lemma 4.4}, there exists $\psi\in L^2_\sigma\left([0,T_0], \mathbb{R}^d\right)$ such that
\[
f_\epsilon^\prime(t)=b_tf_\epsilon(t)+\sigma_t\psi_t  \text{~and~} I(f_\epsilon)=\frac{1}{2}\int_{0}^{T_0}\psi_t^{\mathsf{T}}\sigma_t\psi_tdt.
\]
Let $g_n\in C^2\left([0,T_0], \mathbb{R}^d\right)$ such that $[g_n]$ converges to $[\psi]$ under the distance generated by $\langle \cdot, \cdot\rangle_{L^2_{\sigma, \simeq}}$ as $n\rightarrow+\infty$. For each $n\geq 1$, let $f_n$ be the solution to the ODE
\[
\begin{cases}
&f_n^\prime(t)=b_tf_n(t)+\sigma_tg_n(t) \text{~for~} 0\leq t\leq T_0,\\
&f_n(0)=0,
\end{cases}
\]
then $f_n$ converges to $f_\epsilon$ in $\mathcal{D}\left([0,T_0], \mathbb{R}^d\right)$ and
\[
I(f_n)=\frac{1}{2}\int_0^{T_0} g_n(t)^{\mathsf{T}}\sigma(t)g_n(t)dt=\mathcal{L}_{1,f_n}(g_n)-\frac{1}{2}\mathcal{L}_2(g_n)
\]
by Lemma \ref{lemma 4.4}. Then, there exists $m\geq 1$ such that $f_m\in O$, $I(f_m)\leq I(f_\epsilon)+\epsilon$ and
\begin{equation}\label{equ 4.5}
I(f_m)=\frac{1}{2}\int_0^{T_0} g_m(t)^{\mathsf{T}}\sigma(t)g_m(t)dt=\mathcal{L}_{1,f_m}(g_m)-\frac{1}{2}\mathcal{L}_2(g_m).
\end{equation}
For any $f\in \mathcal{D}\left([0,T_0], \mathbb{R}^d\right)$ and $r>0$, we use $B(f, r)$ to denote the ball concentrated at $f$ with radius $r$.
Since $O$ is open, there exists $\delta(\epsilon)>0$ such that $B\left(f_m,\delta(\epsilon)\right)\subseteq O$. According to the definition of $M_t(f_{n,g_m})$ and the fact that $X_0^n=nx_0=nX_0$,
\begin{align*}
\frac{n}{a_n^2}M_{T_0}(f_{n,g_m})&=\frac{X_{T_0}^n-nX_{T_0}}{a_n}\cdot g_m(T_0)-\frac{X_{0}^n-nX_{0}}{a_n}\cdot g_m(0)\\
&-\int_0^{T_0} \frac{X_s^n-nX_s}{a_n}\cdot g_m^\prime(s) ds-\int_0^{T_0} \sum_{l\in \mathcal{A}}l(\nabla^{\mathsf{T}}F_l)(\xi_{s,l}^n)\frac{X_s^n-nX_s}{a_n}\cdot g_m(s)ds\\
&=\frac{X_{T_0}^n-nX_{T_0}}{a_n}\cdot g_m(T_0)-\int_0^{T_0} \frac{X_s^n-nX_s}{a_n}\cdot g_m^\prime(s) ds\\
&-\int_0^{T_0} \sum_{l\in \mathcal{A}}l(\nabla^{\mathsf{T}}F_l)(\xi_{s,l}^n)\frac{X_s^n-nX_s}{a_n}\cdot g_m(s)ds.
\end{align*}
Then, according to the definitions of $M_{T_0}(f_{n,g_m}), \mathcal{L}_{1,f_m}, \mathcal{L}_2$ and $b_t, \sigma_t$, there exists $\delta_1(\epsilon)\in \left(0, \delta(\epsilon)\right)$ not depending on $n$ such that
\[
\Big|\frac{n}{a_n^2}M_{T_0}(f_{n,g_m})-\mathcal{L}_{1,f_m}(g_m)\Big|<\epsilon\text{\quad and\quad}
\Big|\frac{n}{a_n^2}\langle M(f_{n,g})\rangle_{T_0}-\mathcal{L}_2(g_m)\Big|<\epsilon
\]
when $\vartheta^n\in B(f_m, \delta_1(\epsilon))$ and $\sup_{0\leq t\leq T_0}\left\|\frac{X_t^n}{n}-X_t\right\|\leq \delta_1(\epsilon)$. Then,
by Equation \eqref{equ 4.5} and the expression of $\omega_{T_0}^n(g)$ given in Equation \eqref{equ revision 4.calculationofomega},
\begin{align}\label{equ 4.6}
\left(\omega_{T_0}^n(g)\right)^{-1}&\geq
\exp\left\{\frac{a_n^2}{n}\left(-\mathcal{L}_{1,f_m}(g_m)+\frac{1}{2}(1+o(1))\mathcal{L}_2(g_m)-(2+o(1))\epsilon\right)\right\} \notag\\
&=\exp\left\{\frac{a_n^2}{n}\left(-I(f_m)-(2+o(1))\epsilon+o(1)\mathcal{L}_2(g_m)\right)\right\}
\end{align}
when $\vartheta^n\in B(f_m, \delta_1(\epsilon))$ and $\sup_{0\leq t\leq T_0}\left\|\frac{X_t^n}{n}-X_t\right\|\leq \delta_1(\epsilon)$. We denote by $D^n_{m,\epsilon}$ the event that $\vartheta^n\in B(f_m, \delta_1(\epsilon))$ and $\sup_{0\leq t\leq T_0}\left\|\frac{X_t^n}{n}-X_t\right\|\leq \delta_1(\epsilon)$, then by Lemmas \ref{lemma 4.2 LLNunderPn} and \ref{lemma 4.3 LLNunderGirsanov},
\[
\lim_{n\rightarrow+\infty}P_n^{g_m}\left(D^n_{m, \epsilon}\right)=1.
\]
Therefore, by Equation \eqref{equ 4.6},
\begin{align*}
P(\vartheta^n\in O)&\geq P(D^n_{m, \epsilon})=E_n^{g_m}\Big[\left(\omega_{T_0}^n(g_m)\right)^{-1}1_{\{D^n_{m,\epsilon}\}}\Big] \\
&\geq \exp\big\{\frac{a_n^2}{n}\left(-I(f_m)-(2+o(1))\epsilon+o(1)\mathcal{L}_2(g_m)\right)\big\}(1+o(1)).
\end{align*}
Then,
\[
\liminf_{n\rightarrow+\infty}\frac{n}{a_n^2}\log P(\vartheta^n\in O)\geq -I(f_m)-2\epsilon\geq -I(f_\epsilon)-3\epsilon\geq -\inf_{f\in O}I(f)-4\epsilon.
\]
Since $\epsilon$ is arbitrary, $\liminf_{n\rightarrow+\infty}\frac{n}{a_n^2}\log P(\vartheta^n\in O)\geq -\inf_{f\in O}I(f)$.

\qed

\section{Proof of upper bounds}\label{section five}
In this section we give the proof of the upper bound, where the martingale $\{\omega_t^n(g)\}_{0\leq t\leq T_0}$ introduced in Section \ref{section four} will be utilized. First we show that the upper bound holds for compact sets.

\begin{lemma}\label{lemma 5.1}
For any compact set $\widetilde{K}\subseteq \mathcal{D}\left([0,T_0], \mathbb{R}^d\right)$,
\[
\limsup_{n\rightarrow+\infty}\frac{n}{a_n^2}\log P\left(\vartheta^n\in \widetilde{K}\right)\leq -\inf_{f\in \widetilde{K}}I(f).
\]
\end{lemma}

\proof

According to the definition of $\xi_{t,l}^n, b_t, \sigma_t$, for any $\epsilon>0$ and $g\in C^2\left([0,T_0], \mathbb{R}^d\right)$, there exists $\delta_2>0$ depending on $g$ and $\epsilon$ such that
\[
\left|\int_0^{T_0}g_t^{\mathsf{T}}\left[\sum_{l\in \mathcal{A}}lF_l(\frac{X_t^n}{n})l^{\mathsf{T}}\right]g_tdt-\mathcal{L}_2(g)\right|\leq \epsilon
\]
and
\[
\left|\int_0^{T_0}g_t^{\mathsf{T}}\sum_{l\in \mathcal{A}}l(\nabla^{\mathsf{T}}F_l)(\xi_{t,l}^n)f_tdt-\int_0^{T_0}g_t^{\mathsf{T}}b_tf_tdt\right|\leq \epsilon\sup_{0\leq t\leq T_0}\left\|f_t\right\|
\]
for any $f\in \mathcal{D}\left([0,T_0], \mathbb{R}^d\right)$ when $\sup_{0\leq t\leq T_0}\left\|\frac{X_t^n}{n}-X_t\right\|\leq \delta_2$. Then, conditioned on
$\vartheta^n\in \widetilde{K}$ and $\sup_{0\leq t\leq T_0}\left\|\frac{X_t^n}{n}-X_t\right\|\leq \delta_2$,
\[
\omega_{T_0}^n(g)\geq \exp\left\{\frac{a_n^2}{n}\big[\mathcal{L}_{1,\vartheta^n}(g)-\epsilon\sup_{0\leq t\leq T_0}\left\|\frac{X_t^n-nX_t}{a_n}\right\|
-\frac{1+o(1)}{2}\mathcal{L}_2(g)-\frac{\epsilon(1+o(1))}{2}\big]\right\}
\]
according to the expression of $\omega_{T_0}^n(g)$ given in Equation \eqref{equ revision 4.calculationofomega}. Therefore, by Lemma \ref{lemma 4.1 exponential martingale},
\begin{align*}
1&=E\omega_{T_0}^n(g) \\
&\geq E\Big[\omega_{T_0}^n1_{\{\vartheta^n\in \widetilde{K}\text{~and~}\sup_{0\leq t\leq T_0}\left\|\frac{X_t^n}{n}-X_t\right\|\leq \delta_2\}}\Big] \\
&\geq \exp\left\{\frac{a_n^2}{n}\left[\inf_{f\in \widetilde{K}}\mathcal{L}_{1,f}(g)-\epsilon\sup_{f\in \widetilde{K}}\sup_{0\leq t\leq T_0}\left\|f_t\right\|
-\frac{1+o(1)}{2}\mathcal{L}_2(g)-\frac{\epsilon(1+o(1))}{2}\right]\right\} \\
&\times P\left(\vartheta^n\in \widetilde{K}\text{~and~}\sup_{0\leq t\leq T_0}\left\|\frac{X_t^n}{n}-X_t\right\|\leq \delta_2\right)
\end{align*}
for sufficiently large $n$. As a result,
\begin{align*}
&\limsup_{n\rightarrow+\infty}\frac{n}{a_n^2}\log P\left(\vartheta^n\in \widetilde{K}\text{~and~}\sup_{0\leq t\leq T_0}\left\|\frac{X_t^n}{n}-X_t\right\|\leq \delta_2\right)\\
&\leq -\inf_{f\in \widetilde{K}}\big\{\mathcal{L}_{1,f}(g)-\frac{1}{2}\mathcal{L}_{2}(g)\big\}
+\epsilon\sup_{f\in \widetilde{K}}\sup_{0\leq t\leq T_0}\left\|f_t\right\|+\frac{\epsilon}{2}.
\end{align*}
By Lemma \ref{lemma 3.3} and the fact that $\frac{a_n}{n}\rightarrow 0$,
\[
\limsup_{n\rightarrow+\infty}\frac{n}{a_n^2}\log P\left(\vartheta^n\in \widetilde{K}\text{~and~}\sup_{0\leq t\leq T_0}\left\|\frac{X_t^n}{n}-X_t\right\|\leq \delta_2\right)
=\limsup_{n\rightarrow+\infty}\frac{n}{a_n^2}\log P\left(\vartheta^n\in \widetilde{K}\right).
\]
Hence,
\[
\limsup_{n\rightarrow+\infty}\frac{n}{a_n^2}\log P\left(\vartheta^n\in \widetilde{K}\right)\leq -\inf_{f\in \widetilde{K}}\big\{\mathcal{L}_{1,f}(g)-\frac{1}{2}\mathcal{L}_{2}(g)\big\}
+\epsilon\sup_{f\in \widetilde{K}}\sup_{0\leq t\leq T_0}\left\|f_t\right\|+\frac{\epsilon}{2}.
\]
Since $\epsilon$ and $g$ are arbitrary,
\[
\limsup_{n\rightarrow+\infty}\frac{n}{a_n^2}\log P\left(\vartheta^n\in \widetilde{K}\right)\leq -\left(\sup_{g\in C^2\left([0,T_0], \mathbb{R}^d\right)}\inf_{f\in \widetilde{K}}\big\{\mathcal{L}_{1,f}(g)-\frac{1}{2}\mathcal{L}_{2}(g)\big\}\right).
\]
Note that $\mathcal{L}_{1,f}(g)-\frac{1}{2}\mathcal{L}_{2}(g)$ is convex and continuous of $f$ for fixed $g$ while concave and continuous of $g$ for fixed $f$, then according to the fact that $\widetilde{K}$ is compact and the Minimax Theorem given in \cite{Sion1958},
\[
\sup_{g\in C^2\left([0,T_0], \mathbb{R}^d\right)}\inf_{f\in \widetilde{K}}\big\{\mathcal{L}_{1,f}(g)-\frac{1}{2}\mathcal{L}_{2}(g)\big\}
=\inf_{f\in \widetilde{K}}\sup_{g\in C^2\left([0,T_0], \mathbb{R}^d\right)}\big\{\mathcal{L}_{1,f}(g)-\frac{1}{2}\mathcal{L}_{2}(g)\big\}
=\inf_{f\in \widetilde{K}} I(f)
\]
and the proof is complete.

\qed

To show that the upper bound holds for any closed sets, we need to check that $\{\vartheta^n\}_{n\geq 1}$ is exponential tight. By the main theorem in \cite{Puhalskii1994}, the exponential tightness of $\{\vartheta^n\}_{n\geq 1}$ follows from the following lemma.

\begin{lemma}\label{lemma 5.2} 1)
\begin{equation}\label{equ 5.1}
\limsup_{M\rightarrow+\infty}\limsup_{n\rightarrow+\infty}\frac{n}{a_n^2}\log P\left(\sup_{0\leq t\leq T_0}\left\|\frac{X_t^n-nX_t}{a_n}\right\|>M\right)=-\infty.
\end{equation}
2) For any $\epsilon>0$ and $1\leq i\leq d$,
\[
\limsup_{\delta\rightarrow 0}\limsup_{n\rightarrow+\infty}\frac{n}{a_n^2}\log \sup_{\tau \in \mathcal{T}_0}
P\left(\sup_{0\leq t\leq \delta}\left|e_i^{\mathsf{T}}\cdot\frac{X^n_{\tau+t}-nX_{\tau+t}}{a_n}-e_i^{\mathsf{T}}\cdot\frac{X^n_{\tau}-nX_{\tau}}{a_n}\right|>\epsilon\right)=-\infty,
\]
where $\mathcal{T}_0$ is the set of stopping times of $\{X_t^n\}_{0\leq t\leq T_0}$ with upper bound $T_0$.
\end{lemma}

\proof

For part 1, according to an analysis similar with that leading to Equation \eqref{equ 3.2},
\begin{align}\label{equ 5.2}
&P\left(\sup_{0\leq t\leq T_0}\left\|\frac{X_t^n-nX_t}{a_n}\right\|>M\right) \leq \exp\{-K_3n\}+\sum_{l\in \mathcal{A}}P\left(\sup_{0\leq s\leq nT_1}\big|\widehat{\beta}_l(s)\big|\geq a_nM K_{25}\right) \notag\\
&=\exp\{-K_3n\}+|\mathcal{A}|P\left(\sup_{0\leq s\leq nT_1}\big|\widehat{\beta}(s)\big|\geq a_nM K_{25}\right)
\end{align}
for sufficiently large $n$, where $T_1=K_{11}T_0$, $K_{25}=\frac{e^{-K_9T_0}}{K_{10}|\mathcal{A}|}$ while $\widehat{\beta}(s)=\beta(s)-s$ and $\{\beta(t)\}_{t\geq 0}$ is a Poisson process with rate one. According to an analysis similar with that in the proof of Equation \eqref{equ 3.3},
\[
P\left(\sup_{0\leq s\leq nT_1}\big|\widehat{\beta}(s)\big|\geq a_nM K_{25}\right)
\leq e^{-a_nMK_{25}\theta+nT_1(e^{\theta}-\theta-1)}+e^{-a_nMK_{25}\theta+nT_1(e^{-\theta}+\theta-1)}
\]
for any $\theta>0$. Let $\theta=\frac{a_nMK_{25}}{nT_1}$, then
\[
P\left(\sup_{0\leq s\leq nT_1}\big|\widehat{\beta}(s)\big|\geq a_nM K_{25}\right)\leq \exp\left\{\frac{a_n^2}{n}\left(-\frac{M^2K^2_{25}}{2T_1}+o(1)\right)\right\}.
\]
Therefore, by Equation \eqref{equ 5.2},
\[
\limsup_{n\rightarrow+\infty}\frac{n}{a_n^2}\log P\left(\sup_{0\leq t\leq T_0}\left\|\frac{X_t^n-nX_t}{a_n}\right\|>M\right)\leq -\frac{M^2K^2_{25}}{2T_1}
\]
and then Equation \eqref{equ 5.1} holds.

\quad

For part 2, let $\widehat{D}^n_{M}=\left\{\sup_{0\leq t\leq T_0}\left\|\frac{X_t^n-nX_t}{a_n}\right\|\geq M\right\}\bigcup\left\{\sup_{0\leq t\leq T_0}\left\|X_t^n\right\|\geq nK_2\right\}$, then
\begin{equation}\label{equ 5.3}
\limsup_{M\rightarrow+\infty}\limsup_{n\rightarrow+\infty}\frac{n}{a_n^2}\log P(\widehat{D}_M^n)=-\infty
\end{equation}
by Lemma \ref{lemma 3.2} and Equation \eqref{equ 5.1}. On $\left(\widehat{D}_M^n\right)^c$, by Assumptions (4) and (5), there exists $K_{26}\in (0,+\infty)$ depending on $M$ and $i$ such that
\[
\left|e_i^{\mathsf{T}}\sum_{l\in \mathcal{A}}l(\nabla^{\mathsf{T}}F_l)(\xi^n_{t,l})\frac{X_t^n-nX_t}{a_n}\right|\leq K_{26} \text{~and~}
\left|\left[\sum_{l\in \mathcal{A}}lF_l(\frac{X_t^n}{n})l^{\mathsf{T}}\right]_{ii}\right|\leq K_{26}
\]
for $0\leq t\leq T_0$. Then, for any $\lambda>0$ and sufficiently large $n$,
\begin{align}\label{equ 5.4}
&\left\{\sup_{0\leq t\leq \delta}\left(e_i^{\mathsf{T}}\cdot\frac{X^n_{\tau+t}-nX_{\tau+t}}{a_n}-e_i^{\mathsf{T}}\cdot\frac{X^n_{\tau}-nX_{\tau}}{a_n}\right)>\epsilon\right\}\bigcap \left(\widehat{D}_M^n\right)^c \notag\\
&\subseteq \left\{\sup_{0\leq t\leq \delta}\frac{\omega^n_{t+\tau}(\lambda e_i)}{\omega^n_\tau(\lambda e_i)}
\geq \exp\left[\frac{a_n^2}{n}\left(\lambda\epsilon-\lambda\delta K_{26}-\frac{\lambda^2}{2}(1+o(1))\delta K_{26}\right)\right]\right\}
\end{align}
by taking $g\equiv\lambda e_i$ (and hence $g^\prime=0$) in Equation \eqref{equ revision 4.calculationofomega}. By Lemma \ref{lemma 4.1 exponential martingale}, $\left\{\frac{\omega^n_{t+\tau}(\lambda e_i)}{\omega^n_\tau(\lambda e_i)}\right\}_{0\leq t\leq \delta}$  is a martingale with expectation $1$ for sufficiently large $n$. Then, by Doob's inequality,
\begin{align*}
&P\left(\sup_{0\leq t\leq \delta}\frac{\omega^n_{t+\tau}(\lambda e_i)}{\omega^n_\tau(\lambda e_i)}
\geq \exp\Big[\frac{a_n^2}{n}\left(\lambda\epsilon-\lambda\delta K_{26}-\frac{\lambda^2}{2}(1+o(1))\delta K_{26}\right)\Big]\right)\\
&\leq \exp\left[-\frac{a_n^2}{n}\left(\lambda\epsilon-\lambda\delta K_{26}-\frac{\lambda^2}{2}(1+o(1))\delta K_{26}\right)\right].
\end{align*}
Then, by Equation \eqref{equ 5.4},
\begin{align*}
&\limsup_{n\rightarrow+\infty}\frac{n}{a_n^2}\log \sup_{\tau \in \mathcal{T}_0}
P\left(\sup_{0\leq t\leq \delta}\left(e_i^{\mathsf{T}}\cdot\frac{X^n_{\tau+t}-nX_{\tau+t}}{a_n}-e_i^{\mathsf{T}}\cdot\frac{X^n_{\tau}-nX_{\tau}}{a_n}\right)>\epsilon\right)\\
&\leq \max\left\{-\lambda\epsilon+\lambda\delta K_{26}+\frac{\lambda^2}{2}\delta K_{26}, \text{~}\limsup_{n\rightarrow+\infty}\frac{n}{a_n^2}\log P(\widehat{D}_M^n)\right\}.
\end{align*}
Therefore,
\begin{align*}
&\limsup_{\delta\rightarrow 0}\limsup_{n\rightarrow+\infty}\frac{n}{a_n^2}\log \sup_{\tau \in \mathcal{T}_0}
P\left(\sup_{0\leq t\leq \delta}\left(e_i^{\mathsf{T}}\cdot\frac{X^n_{\tau+t}-nX_{\tau+t}}{a_n}-e_i^{\mathsf{T}}\cdot\frac{X^n_{\tau}-nX_{\tau}}{a_n}\right)>\epsilon\right)\\
&\leq \max\big\{-\lambda\epsilon, \text{~}\limsup_{n\rightarrow+\infty}\frac{n}{a_n^2}\log P(\widehat{D}_M^n)\big\}
\end{align*}
for any $\lambda>0$. Let $\lambda\rightarrow+\infty$, then
\begin{align*}
&\limsup_{\delta\rightarrow 0}\limsup_{n\rightarrow+\infty}\frac{n}{a_n^2}\log \sup_{\tau \in \mathcal{T}_0}
P\left(\sup_{0\leq t\leq \delta}\left(e_i^{\mathsf{T}}\cdot\frac{X^n_{\tau+t}-nX_{\tau+t}}{a_n}-e_i^{\mathsf{T}}\cdot\frac{X^n_{\tau}-nX_{\tau}}{a_n}\right)>\epsilon\right)\\
&\leq \limsup_{n\rightarrow+\infty}\frac{n}{a_n^2}\log P(\widehat{D}_M^n).
\end{align*}
Since $\big\{\frac{\omega^n_{t+\tau}(-\lambda e_i)}{\omega^n_\tau(-\lambda e_i)}\big\}_{0\leq t\leq \delta}$  is also a martingale with expectation $1$ for $\lambda>0$, similar analysis shows that
\begin{align*}
&\limsup_{\delta\rightarrow 0}\limsup_{n\rightarrow+\infty}\frac{n}{a_n^2}\log \sup_{\tau \in \mathcal{T}_0}
P\left(\inf_{0\leq t\leq \delta}\left(e_i^{\mathsf{T}}\cdot\frac{X^n_{\tau+t}-nX_{\tau+t}}{a_n}-e_i^{\mathsf{T}}\cdot\frac{X^n_{\tau}-nX_{\tau}}{a_n}\right)<-\epsilon\right)\\
&\leq \limsup_{n\rightarrow+\infty}\frac{n}{a_n^2}\log P(\widehat{D}_M^n).
\end{align*}
Let $M\rightarrow+\infty$, then part 2 follows from Equation \eqref{equ 5.3}.

\qed

At last we give the proof of the upper bound.

\proof[Proof of the upper bound]

By Lemma \ref{lemma 5.2} and Theorem B on page 47 of \cite{Puhalskii1994}, $\{\vartheta^n\}_{n\geq 1}$ is exponential tight, i.e., for any $m\geq 1$, there exists a compact set $\widetilde{K}_m\subseteq \mathcal{D}\left([0,T_0], \mathbb{R}^d\right)$ such that
\[
\sup_{n\geq 1}\left(P\left(\vartheta^n\notin \widetilde{K}_m\right)^{\frac{n}{a_n^2}}\right)\leq \frac{1}{m}.
\]
For given closed set $C\subseteq \mathcal{D}\left([0,T_0], \mathbb{R}^d\right)$, let $f_m\in C$ such that $\lim_{m\rightarrow+\infty}I(f_m)=\inf_{f\in C}I(f)$.
For each $m\geq 1$, let $\widehat{K}_m$=$\widetilde{K}_m\cup \{f_m\}$, then $\mathrm{\widehat{K}_m}$ is compact and
\begin{equation}\label{equ upperboundproof 1}
\sup_{n\geq 1}\left(P\left(\vartheta^n\notin \widehat{K}_m\right)^{\frac{n}{a_n^2}}\right)\leq \frac{1}{m}
\end{equation}
while $\lim_{m\rightarrow+\infty}\left(\inf_{f\in \widehat{K}_m\cap C}I(f)\right)=\inf_{f\in C}I(f)$. By Lemma \ref{lemma 5.1} and the fact that $\widehat{K}_m\cap C$ is compact,
\[
\limsup_{n\rightarrow+\infty}\frac{n}{a_n^2}\log P\left(\vartheta^n\in C\cap \widehat{K}_m\right)\leq -\inf_{f\in \widehat{K}_m\cap C}I(f).
\]
Then, by Equation \eqref{equ upperboundproof 1},
\[
\limsup_{n\rightarrow+\infty}\frac{n}{a_n^2}\log P\left(\vartheta^n\in C\right)\leq \max\left\{-\inf_{f\in \widehat{K}_m\cap C}I(f), -\log m\right\}
\]
for any $m\geq 1$. Let $m\rightarrow+\infty$, then
\[
\limsup_{n\rightarrow+\infty}\frac{n}{a_n^2}\log P\left(\vartheta^n\in C\right)\leq -\inf_{f\in C}I(f)
\]
follows from the fact that $\lim_{m\rightarrow+\infty}\left(\inf_{f\in \widehat{K}_m\cap C}I(f)\right)=\inf_{f\in C}I(f)$.

\qed

\section{Examples}\label{section six}
In this section we apply our main results in the four examples given in Section \ref{section one}. Throughout this section we assume that $\{a_n\}_{n\geq 1}$ is a positive sequence such that $\lim_{n\rightarrow+\infty}\frac{a_n}{n}=0$ and $\lim_{n\rightarrow+\infty}\frac{a_n^2}{n}=+\infty$.

\textbf{Example 1} \emph{The contact process on the complete graph.} Let $x_0\in (0,1)$ and $X_0^n=nx_0$ for each $n\geq 1$, then
$\{\frac{X_t^n-nX_t}{a_n}\}_{0\leq t\leq T_0}$ follows Theorem \ref{theorem main 2.1 MDP} with
\[
I(f)=\int_0^{T_0}\frac{(f_t^\prime-b_tf_t)^2}{2\sigma_t}dt
\]
for $f$ absolutely continuous, where
\[
X_t=
\begin{cases}
\frac{x_0}{x_0t+1} & \text{~if~}\lambda=1,\\
\frac{(\lambda-1)x_0e^{(\lambda-1)t}}{(\lambda-1)-\lambda x_0+\lambda x_0e^{(\lambda-1)t}} & \text{~if~} \lambda\neq 1,
\end{cases}
\]
$b_t=F_1^\prime(X_t)-F_{-1}^\prime(X_t)=\lambda-2\lambda X_t-1$ and $\sigma_t=F_1(X_t)+F_{-1}(X_t)=X_t(\lambda+1-\lambda X_t)$.

\qed

\textbf{Example 2} \emph{The SIR model on the complete graph.} Let $x_0, y_0$ satisfy $x_0, y_0>0$ while $x_0+y_0<1$ and $X_0^n=\left( nx,  ny_0\right)^{\mathsf{T}}$ for each $n\geq 1$, then $\left\{\frac{X_t^n-nX_t}{a_n}\right\}_{0\leq t\leq T_0}$ follows Theorem \ref{theorem main 2.1 MDP} with
\[
I(f)=\int_0^{T_0}\frac{1}{2}(f_t^\prime-b_tf_t)^{\mathsf{T}}\sigma_t^{-1}(f_t^\prime-b_tf_t)dt
\]
for $f$ absolutely continuous, where $X_t=(S_t, I_t)^{\mathsf{T}}$ satisfies
\[
\begin{cases}
& S_t=x_0e^{-\lambda\phi(t)}, \\
& I_t=-\phi(t)+y_0+x_0(1-e^{-\lambda \phi(t)}),\\
& \phi^\prime(t)=-\phi(t)+y_0+x_0(1-e^{-\lambda \phi(t)}),\\
& \phi(0)=0
\end{cases}
\]
(see the time-change method introduced in Chapter 11 of \cite{Ethier1986}),
\[
b_t=
\begin{pmatrix}
-\lambda I_t & -\lambda S_t\\
\lambda I_t & \lambda S_t-1
\end{pmatrix}
\text{~and~} \sigma_t=
\begin{pmatrix}
\lambda S_tI_t & -\lambda S_tI_t\\
-\lambda S_tI_t & \lambda S_tI_t+I_t
\end{pmatrix}.
\]
Note that it is easy to check that $I_t\geq x_0e^{-t}$ and hence $\sigma_t$ is invertible with
\[
\sigma_t^{-1}=\frac{1}{\lambda S_tI_t^2}\begin{pmatrix}
\lambda S_tI_t+I_t & \lambda S_tI_t\\
\lambda S_tI_t & \lambda S_tI_t
\end{pmatrix}.
\]

\qed

\textbf{Example 3} \emph{Chemical reactions.} Let $x_0, y_0, z_0$ satisfy $x_0, y_0, z_0>0$ while $x_0+y_0+2z_0<1$ and
$X_0^n=\left(nx_0, ny_0, nz_0\right)^{\mathsf{T}}$, then $\left\{\frac{X_t^n-nX_t}{a_n}\right\}_{0\leq t\leq T_0}$ follows Theorem \ref{theorem main 2.1 MDP} with $I(f)$ given by Equation \eqref{equ 2.1 definition of rate function}, where $X_t=(X_t^1, X_t^2, X_t^3)^{\mathsf{T}}$ satisfy
\[
\begin{cases}
& \Big|\frac{X_t^1-c_1}{X_t^1-c_2}\Big|=\Big|\frac{x_0-c_1}{x_0-c_2}\Big|e^{-\lambda(c_1-c_2)t},\\
& X_t^2=X_t^1+y_0-x_0,\\
& X_t^3=x_0+z_0-X_t^1,\\
& c_1, c_2 \text{~are the roots of~}c^2+(y_0-x_0+\frac{\mu}{\lambda})c-\frac{\mu(x_0+z_0)}{\lambda}=0,
\end{cases}
\]
\[
b_t=\begin{pmatrix}
-\lambda X_t^2 & -\lambda X_t^1 & \mu\\
-\lambda X_t^2 & -\lambda X_t^1 & \mu\\
\lambda X_t^2 & \lambda X_t^1 & -\mu
\end{pmatrix}
\]
\text{~and~}
\[
\sigma_t=\begin{pmatrix}
\lambda X_t^1X_t^2+\mu X_t^3 & \lambda X_t^1X_t^2+\mu X_t^3 & -\left(\lambda X_t^1X_t^2+\mu X_t^3\right)\\
\lambda X_t^1X_t^2+\mu X_t^3 & \lambda X_t^1X_t^2+\mu X_t^3 & -\left(\lambda X_t^1X_t^2+\mu X_t^3\right)\\
-\left(\lambda X_t^1X_t^2+\mu X_t^3\right) & -\left(\lambda X_t^1X_t^2+\mu X_t^3\right) &  \lambda X_t^1X_t^2+\mu X_t^3
\end{pmatrix}.
\]
It is easy to check that $I(f)<+\infty$ implies $f_t=(f_1(t), f_1(t), -f_1(t))^{\mathsf{T}}$ for some absolutely continuous $f_1(t):[0,T_0]\rightarrow \mathbb{R}$ and then
\[
I(f)=\frac{1}{2}\int_0^{T_0}\psi_t^{\mathsf{T}}\sigma_t\psi_tdt
\]
with $\psi(t)=(\psi_1(t), \psi_1(t), -\psi_1(t))^{\mathsf{T}}$ and $\psi_1(t)=\frac{f_1^\prime(t)+(\lambda
X_t^1+\lambda X_t^2+\mu)f_1(t)}{3(\lambda X_t^1X_t^2+\mu X_t^3)}$ by Lemma \ref{lemma 4.4}.

\qed

\textbf{Example 4} \emph{Yule process with rate $\lambda$.} Let $x_0>0$ and $X_0^n=nx_0$ for each $n\geq 1$, then $\left\{\frac{X_t^n-nX_t}{a_n}\right\}_{0\leq t\leq T_0}$ follows Theorem \ref{theorem main 2.1 MDP} with
\[
I(f)=\int_0^{T_0}\frac{(f_t^\prime-b_tf_t)^2}{2\sigma_t}dt=\int_0^{T_0} \frac{(f_t^\prime-\lambda f_t)^2}{2\lambda x_0e^{\lambda t}}dt
\]
for $f$ absolutely continuous, where $X_t=x_0e^{\lambda t}$, $b_t\equiv \lambda$ and $\sigma_t=\lambda X_t=\lambda x_0e^{\lambda t}$.

\qed

\quad

\textbf{Acknowledgments.} The author is grateful to Dr. Linjie Zhao for useful suggestions. The author is grateful to the financial
support from the National Natural Science Foundation of China with
grant number 11501542.

{}
\end{document}